\documentclass{nm}
\renewcommand\thisnumber{x}
\renewcommand\thisyear {200x}

\renewcommand\thisvolume{xx}
\usepackage{charter}
\usepackage[charter]{mathdesign}

\usepackage{bbm}
\usepackage{booktabs}
\usepackage{amsmath}

\newcommand{\Partial}[3]{\ensuremath{\frac{\partial^{#1}{#2}}{\partial{#3}^{#1}}}}
\renewcommand{\Partial}[3][ ]{\ensuremath{\partial^{#1}_{#3} {#2}}}

\newcommand {\bfvec}[1]{\ensuremath{{\mathbf #1}}}
\newcommand{\vecx}{\ensuremath {\bfvec x}}

\newcommand{\grad}{\ensuremath{ \nabla }}

\newcommand{\Div}[1]{\ensuremath{ \nabla \cdot \left( #1 \right)}}

\newcommand{\clos}[1]{\ensuremath{\overline{#1}}}
\newcommand{\dx}{\ensuremath{ \,{\rm d}x  } }
\newcommand{\dS}{\ensuremath{{\rm d}S}}

\newcommand{\Gammat}{\ensuremath{{{\mathcal G}^t}}}
\newcommand{\Gammanul}{\ensuremath{{{\mathcal G}^0}}}

\newcommand{\R}[1]{\ensuremath{\mathbbm{R}^{#1}}}

\newcommand{\tldvecu}{\ensuremath{\tilde {\bfvec u}}}
\newcommand{\hatvecu}{\ensuremath{\hat {\bfvec u}}}

\newcommand{\Omegain}{\ensuremath{\Omega_{in}}}
\newcommand{\Omegaout}{\ensuremath{\Omega_{out}}}

\begin{document}

\markboth{V. Klement, T. Oberhuber, D. \v{S}ev\v{c}ovi\v{c}}{Application of the level-set model with constraints in image segmentation}
\title{Application of the level-set model with constraints in image segmentation}

\author[Author(s)]{Vladim\'{i}r Klement\affil{1}, Tom\'{a}\v{s} Oberhuber\affil{1} and Daniel \v{S}ev\v{c}ovi\v{c}\affil{2}\comma\corrauth}

\address{\affilnum{1}\ Department of Mathematics, Faculty of Nuclear Sciences and Physical Engineering, Czech Technical University in Prague, Trojanova 13, Praha 2, 120 00, Czech Republic\\
\affilnum{2}\ Department of Applied Mathematics and Statistics, Comenius University, 842 48 Bratislava, Slovakia}

\emails{{\tt tomas.oberhuber@fjfi.cvut.cz (T. Oberhuber), vladimir.klement@fjfi.cvut.cz (V. Klement), sevcovic@fmph.uniba.sk (D. \v{S}ev\v{c}ovi\v{c})} }

\begin{abstract}
We propose and analyze a constrained level-set method for semi-automatic image segmentation. Our level-set model with constraints on the level-set function enables us to specify which parts of the image lie inside respectively outside the segmented objects. Such a-priori information can be expressed in terms of upper and lower constraints prescribed for the level-set function. Constraints have the same conceptual meaning as initial seeds of the popular graph-cuts based methods for image segmentation. A numerical approximation scheme is based on the complementary-finite volumes method  combined with the Projected successive over-relaxation method adopted for solving constrained linear complementarity problems. The advantage of the constrained level-set method is demonstrated on several artificial images as well as on cardiac MRI data.
\end{abstract}

\keywords{linear complementarity, image processing, segmentation, level-set method, projected successive over-relaxation method}

\ams{90C33, 35K55, 35K52, 53C44, 74S10, 74G15}

\maketitle

\section{Introduction}

{\it The level-set methods} for the image segmentation have been studied and applied during the last two decades. The level-set method applied in the image segmentation is typically an iterative method. The segmentation starts with an initial curve $\Gammanul$ representing an initial guess for the segmented object and it is evolved in the normal direction towards the segmented object by means of a suitable geometric law taking into account the orientation of the segmented object and also the curvature of evolved curves. Loosely speaking, the better the initial guess is, the better and faster the segmentation process is. This is profitable for processing of time sequences where the final segmentation of one frame may serve as the initial guess for the next frame. We refer the reader to a wide range of literature on this topic e.g. Caselles et al. \cite{CasellesCatteCollDibos-1993}, Handlovi\v{c}ov\'{a} et al. \cite{HandlovicovaMikulaSgallari-2003}, Osher, Paragios \cite{OsherParagios-2003} or Sethian \cite{Set_LSMEIIGFMCV&MS} and references therein. In comparison to parametric models studied by Bene\v{s} et al. \cite{BenesKimuraPausSevcovicTsujikawaYazaki-2008} and Kass et al. \cite{KassWitkinTerzopoulos-1987} the level-set methods can handle topological changes and therefore one initial curve can split and segment more separate objects. The level-set method is still subject of very active research. In \cite{BenesMaca-2012}, time sequences of 2D MRI slices are segmented as 3D data by the level-set method. It ensures smooth segmentation of adjacent slices. The multilayer segmentation level-set method for segmentation of images with nested structures is presented in \cite{ChungVese-2009}. Combination of the level-set methods with statistical approaches is subject of the review paper \cite{CremersRoussonDeriche-2007}. Review of deformable contour models in medical image segmentation can be found in \cite{HePengEverdingWangHanWeissWee-2008}.

Among different segmentation methods there are {\it the graph-cuts methods} (see e.g. Boykov et al. \cite{BoykovFunka-Lea-2006,BoykovKolmogorov-2004}, Gurholt and Tai \cite{Gurholt-2009}, Louck\'y and Oberhuber \cite{LouckyOberhuber-2012}) which are based on the graph theory and algorithms for finding {\it minimal cuts} and {\it the maximal flow} respectively. These algorithms are not iterative and they do not require initial curves. Instead of it, they need initial seeds - one or more points or lines in the interior and exterior of the segmented object.

Each segmentation algorithm requires some description of the object of our interest. The object is described usually in some of the following ways:
{\it Edges} -- it is often used information since many objects in the real world have clearly visible edges. In the level-set methods, the Perrona-Malik function serves as an edge detector.
{\it Color or texture pattern} -- real objects usually have uniform or homogeneous surface. Therefore areas of the same color or texture pattern belongs very likely to objects of the same type.
{\it Shape} -- another criterion might be segmentation of objects with prescribed shape. The object shape can be given by insertions of an appropriate anisotropy \cite{Oberhuber-2009}, the shape-learning methods or by minimizing the elastic energy of the segmentation curve \cite{DroskeRumpf-2004}.
{\it Location} -- it is expressed by the initial condition. Proper setting of the initial curve for the level-set segmentation may help to specify what object we aim to segment especially if there are more similar objects. Note however, that the initial curve of the level-set method is only an initial step for the segmentation algorithm and the final segmentation may differ from the initial curve significantly.
{\it Skeleton} -- the initial seeds in the graph-cuts method differ from the initial curve in one important fact. What is marked by the initial seed as an interior of the segmented object will remain interior even in the final segmentation and vice versa for the exterior.

From this point of view, we can understand the initial seeds in the graph-cuts method as a {\it hard segmentation constraint} while the initial curve in the level-set method as a {\it soft segmentation constraint}. In this article, we show how to incorporate local a-priori information similar to the initial seeds used in the graph-cuts method to the level-set method. We propose a new constrained level-set method which can be applied to the image segmentation problems. For better understanding of our method, we will compare it with the classical level-set methods (c.f.  \cite{CasellesCatteCollDibos-1993}) with no surface terms extracting the information about the object color or texture. 

The constrained level set method allows an expert to prescribe an a-priori information by marking parts which are surely inside or outside the segmented region. As an example of application of the constrained level set method we chose cardiac medical images shown in Figures~\ref{fig:cardiac-mri-1},\ref{fig:cardiac-mri-2}. Applying our level-set method with constraints, we may mark septum (red lines in Figure~\ref{fig:cardiac-mri-1} a)) as "must stay outside the segmented region" and we can obtain correct result. On the other hand, the unconstrained level-set method Figure~\ref{fig:cardiac-mri-1} b) was not capable to segment the left and the right ventricle separately. An advantage of the proposed method, compared to the graph-cuts methods, is that it allows to incorporate anisotropies \cite{Oberhuber-2009} or other energies to minimize like the elastic one \cite{DroskeRumpf-2004}. Notice that our aim is to compare segmentation results obtained by the constrained and unconstrained level set method as well as the graph-cuts method. It should be emphasized that for such a specific cardiac segmentation problem there are other more sophisticated and fully unsupervised methods utilizing specific information about the image. There are also level set implementation of the region based methods (see e.g. \cite{ParagiosDeriche-2000}, \cite{Jiang-2012} and references therein). Nevertheless, we do not present comparison to those specific methods. 

A numerical approximation scheme is based on the complementary-finite volumes method developed by Handlovi\v{c}ov\'a et al. in  \cite{HandlovicovaMikulaSgallari-2003} combined with the projected successive over-relaxation method for solving constrained problems proposed by Mangasarian in \cite{Mangasarian-1977} and Elliott, Ockendon in \cite{ElliottOckendon-1982}. The advantage of the constrained level-set method is demonstrated by means of several artificial images. 

The paper is organized as follows. In Section~\ref{sect:level-set-segmentation} we recall general level-set method for the image segmentation together with the numerical scheme and successive over-relaxation (SOR) method. In Section~\ref{sect:constrained-level-set}, we present our level-set model with constraints together with a efficient numerical scheme. As a solver for the linear complementarity problem with range bounds we adopt the Projected SOR method. Comparison with the common level-set method and contributions of the constrained level-set method are demonstrated in Section~\ref{sect:results}.

\section{The level-set method for the image segmentation}
\label{sect:level-set-segmentation}

\subsection{Time-space continuous framework of the level-set method}

We consider a given image which is represented by the greyscale image function $I_0:\Omega\to [0,1]$ defined on a two dimensional rectangle $\Omega \equiv [K_1,L_1]\times[K_2,L_2]$. The idea how to segment an object in the image is to start from a closed, embedded and smooth initial curve $\Gammanul$ approximating the shape of the object and let it evolve towards the exact boundary of the object. To this end, we construct a family of evolving curves $\Gammat$ with the property that $\Gammat$ converges to the boundary of a segmented object as $t$ goes to infinity. There are many ways how to construct such a flow of planar curves. Among them we will focus our attention to the flow of curves proposed in the active contour model (c.f. Caselles et al.  \cite{CasellesKimmelSapiro-1995}, Kichenassamy et al. \cite{KKOTY2}). A problem of finding a boundary of an object in the image can be reformulated as a problem of construction of planar curves on which the gradient $\nabla I_0$ of the image intensity function $I_0$ is large.

Assuming $\Gammat$ is a $C^1$ smooth curve we can evaluate the unit tangent vector $\bfvec T(\bfvec x)$ and outer unit normal vector $\bfvec N(\bfvec x)$. Each point $\bfvec x \in \Gammat$ is evolved in the normal direction with the normal velocity $V=V(\bfvec x,t)$, i. e.
 \[
 \partial_t \bfvec x = V(\bfvec x,t) . N(\bfvec x,t).
 \]
Although the velocity vector $\partial_t \bfvec x$ can be decomposed into its tangential and normal parts, it should be noted that only the motion in the normal direction has impact on the shape of the closed curve $\Gammat$. 

Following the active contour model (c.f. \cite{CasellesKimmelSapiro-1995,KKOTY2}), Mikula and \v{S}ev\v{c}ovi\v{c}, in \cite{MikulaSevcovic2004} Mikula and \v{S}ev\v{c}ovi\v{c} considered a generalized form of the normal velocity:
\begin{equation}
\label{eqn:normal-velocity}
V\left(\bfvec x, t\right) = g^0\left(\bfvec x\right) H\left(\bfvec x, t\right) + \grad g^0\left(\bfvec x\right) \cdot \bfvec N\left(\bfvec x, t\right),
\end{equation}
where $H(\bfvec x,t)$ is the curvature of $\Gammat$. It is known that it has a smoothing effect on the segmented curve and the curvature driven flow is the gradient flow for the total length of a curve (c.f. \cite{MikulaSevcovic2004}). Next, $g^0 = g\left(\left| G_\sigma \ast \grad I_0\right| \right)$ where $g$ is a smooth edge detector function $g:[0,\infty)\to (0,\infty)$ such that  $g^\prime<0, g(0)=1, g(+\infty)=0$ and $g^\prime(s)\le C g(s),|g''(s)|\le C, s>0$, for some constant $C>0$. A typical example is the function $g(s)=1/(1+\lambda s^2)$ where $\lambda>0$ is a contrast parameter. Notice that, for a given smooth intensity function $I_0$, the vector field  $\vec W(x) = - \nabla g^0(x)$  has an important geometric property as it points towards edges in the image where the norm of the gradient $\nabla I_0$ is large (c.f. \cite{MikulaSevcovic2004}). Notice that a possible lack of smoothness of $I_0$ (e.g. due to a noise) can be overcome by taking the convolution of $I_0$ with a smooth Gaussian mollifier $G_\sigma$ with the variance $\sigma^2>0$ (see \cite{MikulaSevcovic2004,MikulaSarti-2007}). The term $\grad g^0\left(\bfvec x\right) \cdot \bfvec N\left(\bfvec x, t\right)$ pushes the evolved curve $\Gammat$ towards the edge of the image $I_0$ (c.f. \cite{KassWitkinTerzopoulos-1987,MikulaSevcovic2004}). The effect of the curvature term $H\left(\bfvec x, t\right)$ consists in smoothing the segmented curve by means of minimization of its total length. This property makes the segmentation model robust for application even in the case of a noisy image. Notice that the term $g^0$ slows down the normal velocity in the vicinity of edges of $I_0$  (c.f. \cite{MikulaSevcovic2004}). 

In the level-set method, $\Gammat$ is given implicitly as
%
%\begin{equation}                zruseno cislovani Tomas 3.5.
%\label{def:level-set-gamma}
$$
\Gammat \equiv \left\{ \bfvec x \in \Omega \mid u\left(\bfvec x,t\right) = 0 \right\},
$$
%\end{equation}
%
where $u$ is a real valued smooth function defined on $\Omega$ such that $u\left(\bfvec x \right)<0$ for all $\bfvec x$ belonging to the interior of $\Gammat$ and $u\left(\bfvec x\right)>0$ for all $\bfvec x$ belonging to the exterior of a Jordan curve $\Gammat$. Following derivation from \cite{MikulaSarti-2007}, the level-set formulation of (\ref{eqn:normal-velocity}) can be stated in terms of a solution $u$ to the initial-boundary value problem
\begin{eqnarray}
\label{eq:level-set-segmentation}
&u_t = Q \Div{g^0 \frac{\grad u}{Q}}\ &{\rm in}\ \Omega \times (0, T], \\
&\Partial{u}{\nu} = 0\ &{\rm on}\ \partial \Omega, \\
\label{eq:level-set-segmentation-3}
&u \mid_{t=0} = u_{ini}\ &{\rm in}\ \Omega,
\end{eqnarray}
where $u_{ini}$ is the initial level-set function corresponding to the initial curve $\Gammanul$, $\Partial{u}{\nu} = \grad u \cdot \nu$ and $\nu$ is the outer normal unit vector of the boundary $\partial \Omega$ of a computational domain $\Omega$. In the level-set framework, the quantity $Q$ should be equal to $\left| \grad u \right|$. However, for practical purposes, it is regularized by means of the Tichonov regularization, i.e. 
\begin{equation}
\label{def:Q}
Q=\sqrt{\epsilon^2 + \left| \grad u \right|^2},
\end{equation}
where $0<\epsilon\ll 1$ is a small regularizing parameter. 

\subsection{Time-space discretized framework of the level-set method}

\begin{figure}
\center{
\includegraphics[width=8cm]{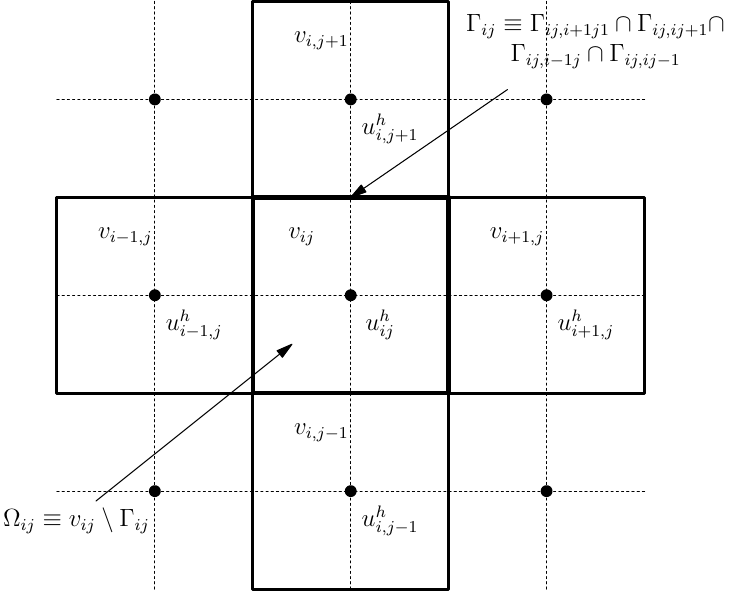}
}
\caption{A description of a finite volume $v_{ij}$ defined on the dual mesh. Here $\Omega_{ij}$ is its interior, $\Gamma_{ij}$ its boundary consisting of linear segments $\Gamma_{ij,i\pm1j}, \Gamma_{ij,ij\pm1}$.}
\label{fig:fin-volume}
\end{figure}

We discretize the initial-boundary value problem  (\ref{eq:level-set-segmentation})--(\ref{eq:level-set-segmentation-3}) by means of the method of complementary finite volumes developed by Handlovi\v{c}ov\'a et al. \cite{HandlovicovaMikulaSgallari-2003} in the context of a class of level-set equations arising in the image processing. Let $\tau$ be a time step for time discretization. Let $h=(h_1, h_2)$ be spatial discretization steps such that $h_i = \frac{L_i-K_i}{N_i}$ for some $N_i\in \mathbbm{N}^+, i=1,2$. We define a numerical grid 
\[
M_h = \left\{ (ih_1, jh_2) \mid i = 0, \cdots, N_1, j = 0, \cdots, N_2 \right\}.
\]
For a function $u \in C\left(\clos{\Omega}\times(0,T];\R{}\right)$ we define its piecewise constant approximation on $M_h$ at the time $k \tau$ as a grid function defined by $u^k_{ij} = u\left(ih_1,jh_2,k\tau\right)$. We furthermore introduce a dual mesh $V_h$ defined as
\begin{eqnarray*}
V_h &\equiv& \left\{ v_{ij} =  \left[ \left(i-\frac{1}{2}\right) h_1, \left(i+\frac{1}{2}\right) h_1 \right] \times
                             \left[ \left(j-\frac{1}{2}\right) h_2, \left(j+\frac{1}{2}\right) h_2 \right] \mid  \right. \nonumber \\
                             &&  i = 1, \cdots, N_1-1, j = 1, \cdots, N_2-1 \bigg\}.
\label{dual-fvm-mesh}
\end{eqnarray*}
For $0 < i < N_1$, $0 < j < N_2$, $i$ and $j$ fixed, we consider a finite volume $v_{ij}$ of the dual mesh $V_h$. We denote its interior by $\Omega_{ij}$, its boundary by  $\Gamma_{ij}$ and let $\mu\left(v_{ij}\right)$ be the area of $\Omega_{ij}$ (see Fig.~\ref{fig:fin-volume}).  We also denote the set of all neighboring volumes (having one common edge) of a finite volume $v_{ij}$ by $\mathcal{N}_{ij}$. For all inner finite volumes $v_{ij}$ of the dual mesh $V_h$, the boundary $\Gamma_{ij}$ consists of four linear segments. We denote them as $\Gamma_{ij,\bar{ij}}$ for $\bar{ij} \in \mathcal{N}_{ij}$. It means that $\Gamma_{ij,\bar{ij}}$ is a boundary between the finite volumes $v_{ij}$ and $v_{\bar{ij}}$. Let $l_{ij,\bar{ij}}$ be the length of this part of $\Gamma_{ij}$.
Dividing equation (\ref{eq:level-set-segmentation}) by the term $Q$, integrating it over the interior $\Omega_{ij}$ of a finite volume $v_{ij}$ we obtain the equation:
%
%\begin{equation}  zruseno cislovani - 3.5. Tomas
$$
\int_{\Omega_{ij}}\frac{u_t}{Q} \dx = \int_{\Omega_{ij}}\Div{g^0 \frac{\grad u}{Q}} \dx.
$$
%\end{equation}
%
Applying the Euler backward difference in time discretization of $u_t$ we end up with the following time semi-discretization of (\ref{eq:level-set-segmentation}):
%
%\begin{equation}  zruseno cislovani - 3.5. Tomas
$$
\int_{\Omega_{ij}}\frac{1}{Q^{k-1}} \frac{u^{k}-u^{k-1}}{\tau}\dx = \int_{\Omega_{ij}}\Div{g^0 \frac{\grad u^k}{Q^{k-1}}}\dx.
$$
%\end{equation}
%
Next, after long but straightforward calculations and evaluations described in Appendix we are in a position to formulate full time-space discretization of the level set equation (\ref{eq:level-set-segmentation}). It is a semi-implicit scheme in which the nonlinear terms are treated explicitly. It can be rewritten in a form of the following system of linear equations:
\begin{equation}
A_{ij}^{k} u_{ij}^{k} +
A_{i+1j}^{k} u_{i+1j}^{k} +
A_{ij+1}^{k} u_{ij+1}^{k} +
A_{i-1j}^{k} u_{i-1j}^{k} +
A_{ij-1}^{k} u_{ij-1}^{k} = u_{ij}^{k-1},
\label{eq:algebraic-segmentation}
\end{equation}
for $i=1,\cdots, N_1-1$ and $j=1,\cdots, N_2-1$, and 
\begin{eqnarray*}
A^k_{0j} u^k_{0j} + A^k_{1j} u^k_{1j} &=& 0\ {\rm for}\ j=0, \cdots, N_2, \\
A^k_{N_1j} u^k_{N_1j} + A^k_{N_1-1j} u^k_{N_1-1j} &=& 0\ {\rm for}\ j=0, \cdots, N_2, \\
A^k_{i0} u^k_{i0} + A^k_{i1} u^k_{i1} &=& 0\ {\rm for}\ i=0, \cdots, N_1, \\
A^k_{iN_2j} u^k_{iN_2} + A^k_{iN_2-1} u^k_{iN_2-1} &=& 0\ {\rm for}\ i=0, \cdots, N_1,
\end{eqnarray*}
where the terms $A^k_{ij}$ are derived in Appendix (see  (\ref{def:A-1})--(\ref{def:A-6}) ).

At each time level $k\tau$ we can represent a solution $u^k_{ij}$ by a stacked vector 
\[
\tldvecu\equiv\tldvecu^k=\left(u^k_{00}, \cdots, u^k_{0N_2}, u^k_{10}, \cdots, u^k_{1N_2}, u^k_{20}, \cdots, u^k_{N_10}, \cdots, u^k_{N_1N_2}\right)^T,
\]
by mapping the node $(i,j)$ of the two dimensional spatial domain to the one-dimensional vector, i.e. $I=I(i,j)=j \cdot N_1 + i$ for $i=0,\cdots,N_1$ and $j=0,\cdots,N_2$ and setting $\tldvecu^k_I = u^k_{ij}$. With this notation, the system of linear equations (\ref{eq:algebraic-segmentation}) can be then rewritten in a matrix form
\begin{equation}
\label{eq:linear-problem}
\mathbbm{A} \tldvecu = \bfvec b
\end{equation}
for the solution vector $\tldvecu\equiv\tldvecu^k$. The dimension of the square matrix $\mathbbm{A}$ as well as of the vectors $\tldvecu$ and $\bfvec b$ is $(N_1+1)(N_2+1)$. The $I=I(i,j)$-th row of the matrix $\mathbbm{A}=(a_{IJ})$ contains only five nonzero elements $A^k_{ij}, A^k_{i\pm1 j}, A^k_{ij\pm1 }$. Now it follows from (\ref{def:A-1})--(\ref{def:A-6}) that $A^k_{i\pm1 j}, A^k_{ij\pm1} <0$ and the diagonal term $a_{II} = A_{ij}^{k}   = 1 - ( A_{i+1j}^{k} + A_{ij+1}^{k} + A_{i-1j}^{k} + A_{ij-1}^{k})>0$. Hence the matrix $\mathbbm{A}$ is a sparse diagonally dominant $M$-matrix. We can solve the problem (\ref{eq:linear-problem}) by means of the {\it Successive over-relaxation method} (SOR) \cite{Quarteroni-2000}. 

\section{The constrained level-set method and its numerical approximation}
\label{sect:constrained-level-set}

\subsection{Time-space continuous framework of the constrained level-set method}

In this section we introduce our constrained level-set method for image segmentation. We suppose that there are two disjoint subdomains $\Omegain$ and $\Omegaout$ of the domain $\Omega$ such that $\Omegain$ is a subset of the interior of the segmented object and $\Omegaout$ lies outside the segmented object. Furthermore, we suppose that there are two prescribed functions $v,w \in C\left(\Omega\right)$ with the property such that $w < v$ in $\Omega$ and $v<0$ in $\Omegain$ and $v>0$ in $\Omega \setminus \Omegain$ and $w>0$ in $\Omegaout$ and $w<0$ in $\Omega \setminus \Omegaout$. Notice that any function $u$ fulfilling $w \leq u \leq v$ in $\Omega$ must be negative in $\Omegain$ and positive in $\Omegaout$. Its zero level-set contains the set $\Omegain$ in its interior and $\Omegaout$ in its exterior. 

Our purpose is to construct a solution $u=u(x,t)$ such that it satisfies the level set equation (\ref{eq:level-set-segmentation}) in the open region where $w(x) < u(x,t) < v(x)$. Moreover, we require that $w(x) \le  u(x,t) \le v(x)$ for all $x\in\Omega, t\in(0,T]$. 
The reason why to prescribe range bounds on the level-set function $u$ is to keep the set  $\Omegain$ inside and $\Omegaout$ outside the zero level set ${\mathcal G}^t=\{\vecx, u(\vecx,t)=0\}$ approaching the boundary of a segmented object when $t\to\infty$. To this end we consider the following partial differential inequality problem:
\begin{eqnarray}
u_t &=& Q \Div{g^0 \frac{\grad u}{Q}}\ 
\hbox{for}\ (x,t)\in\Omega \times (0, T],\ \hbox{such that}\ w(x) < u(x,t) < v(x), 
\nonumber
\\
\label{eq:constrained-level-set-segmentation}
\\
u_t &\ge&  Q \Div{g^0 \frac{\grad u}{Q}}\ 
\hbox{for}\ (x,t)\in\Omega \times (0, T],\ \hbox{such that}\ u(x,t)=w(x), 
\\
u_t &\le&  Q \Div{g^0 \frac{\grad u}{Q}}\ 
\hbox{for}\ (x,t)\in\Omega \times (0, T],\ \hbox{such that}\ u(x,t) = v(x), 
\nonumber
\end{eqnarray}
\begin{eqnarray}
\Partial{u}{\nu} &=& 0\ {\rm at}\ \partial \Omega, \\
\label{eq:constrained-level-set-segmentation-3}
u \mid_{t=0} &=& u_{ini}\ {\rm in}\ \Omega.
\end{eqnarray}
The precise mathematical formulation of (\ref{eq:constrained-level-set-segmentation}) can be stated in terms of the following variational inequality problem: given barrier functions $v, w \in W^{1,2}\left(\Omega\right)$, $w(x) < v(x)$ for a.e. $x \in \Omega$, find a solution $u\in {\mathcal K}\subset {\mathcal X}$ where ${\mathcal X} = W^{1,2}((0,T): L^2(\Omega)) \cap L^2((0,T): W^{1,2}(\Omega))$ such that  
\[
\langle{\mathcal A}(u), \phi-u \rangle \ge 0, \quad \hbox{for each}\ \phi \in {\mathcal K},
\]
where ${\mathcal K}$ is the convex closed cone (c.f. Brezis \cite[Chapter 2]{Brezis1972}):
\begin{equation}
{\mathcal K} = \{ u \in {\mathcal X}, \ w(x) \leq u(x, t) \leq v(x), \ \hbox{for a.e.}\ x\in \Omega, t\in(0,T)\}.
\label{ds:A1}
\end{equation}
Here $L^p(\Omega), 1\le p<\infty$ is the usual Lebesgue space 
 consisting of all measurable functions on $\Omega$ such that $\Vert 
 u\Vert_p = (\int_\Omega|f(x)|^p dx)^{1/p} < \infty$. Furthermore, 
 $W^{1,2}(\Omega)$ stands for the Sobolev space consisting of all 
 functions having finite Sobolev norm $\Vert u\Vert_{1,2} = \Vert 
 u\Vert_2 + \Vert Du\Vert_2$ where $Du$ is the gradient of $u$ in the 
 sense of distributions. Finally,  $W^{-1,2}(\Omega)$ is the dual space 
 to  $W^{1,2}(\Omega)$.
The operator  ${\mathcal A}: {\mathcal X} \to  L^{2}((0,T): W^{-1,2}(\Omega))$ is defined by
\[
{\mathcal A}(u) = \frac{u_t}{Q} -   \Div{g^0 \frac{\grad u}{Q}}
\]
and $\langle ., . \rangle$ is the inner product in $L^{2}((0,T): L^{2}(\Omega))$, i.e. 
\[
\langle{\mathcal A}(u), \phi-u \rangle = \int_0^T \int_\Omega 
\left( u_t (\phi-u) +   g^0 \grad u \cdot \grad (\phi-u) \right)\frac{\dx}{Q} dt,
\]
where $Q=Q\left(u\right)=\sqrt{\epsilon^2+\left|\grad u\right|^2}$, $g^0=g^0\left(u\right)$.

\subsection{Time-space discretized framework of the constrained level-set method}

The discretization of (\ref{eq:constrained-level-set-segmentation})--(\ref{eq:constrained-level-set-segmentation-3}) follows exactly from the discretization of the level set equation (\ref{eq:level-set-segmentation})--(\ref{eq:level-set-segmentation-3}) when taking into account the range bound constraints $w(x) < u(x, t) < v(x)$ for a.e. $x\in\Omega, t\in(0,T)$. 
At each time step we have to construct a solution $\tldvecu$ to the following linear complementarity problem: 
\begin{eqnarray}
(\mathbbm{A} \tldvecu)_I &=& \bfvec b_I \ \hbox{for $I$ such that} \ w_I < \tldvecu_I < v_I,
\nonumber
\\
(\mathbbm{A} \tldvecu)_I &\ge & \bfvec b_I \ \hbox{for $I$ such that} \ \tldvecu_I =w_I,
\label{eq:linear-problem-inequality}
\\
(\mathbbm{A} \tldvecu)_I &\le & \bfvec b_I \  \hbox{for $I$ such that}\ \tldvecu_I =v_I,
\nonumber
\end{eqnarray}
for $I=1, \cdots, (N_1+1)(N_2+1)$.

In order to solve the linear complementarity problem (\ref{eq:linear-problem-inequality}) we make use of the so-called Projected SOR method (PSOR) \cite{Mangasarian-1977,ElliottOckendon-1982} adopted for our problem. 

For each index $I=1, \cdots, (N_1+1)(N_2+1)$, we repeat the following up-dating of the vector $\tldvecu^{(p)}$: 
\begin{eqnarray}
\hatvecu^{(p+1)}_I &=& 
\left(1 - \omega \right)\tldvecu_I^{(p)} + \frac{\omega}{a_{II}}\left( 
\bfvec b_I - \sum_{J<I} a_{IJ} \tldvecu_J^{(p+1)} - \sum_{J>I} a_{IJ} \tldvecu_J^{(p)} \right), 
\label{def:psor-sor}
\\
\tldvecu^{(p+1)}_I &=& 
\min\{
     \max\{\hatvecu^{(p+1)}_I,  w_I\}, v_I
    \}, 
\label{def:psor}
\end{eqnarray}
until the prescribed tolerance level for the difference 
$\Vert \tldvecu^{(p+1)} - \tldvecu^{(p)}\Vert < tol$ is achieved. Equation (\ref{def:psor-sor}) corresponds to a usual SOR method. Equation (\ref{def:psor}) ensures that the prescribed constraints are satisfied. Recall that the matrix $\mathbbm{A}$ is a positively diagonally dominant $M$-matrix and so $a_{II}>0$ for each $I$. Assume that $\tldvecu^{(p)} \to \tldvecu$ as $p\to\infty$. Hence $\hatvecu^{(p+1)} \to \left(1 - \omega \right)\tldvecu + \omega \mathbbm{D}^{-1}\left( \bfvec b - (\mathbbm{L}+\mathbbm{U})\tldvecu  \right)$. 
It means 
\begin{equation}
\tldvecu_I =
\min\{
     \max\{
[ \left(1 - \omega \right)\tldvecu + \omega \mathbbm{D}^{-1}\left( 
\bfvec b - (\mathbbm{L}+\mathbbm{U})\tldvecu  \right)]_I,  w_I\}, v_I
    \}.
\label{eqA}
\end{equation}
Clearly, 
\[
w_I \le \tldvecu_I \le v_I
\]
for each index $I$. Here $\mathbbm{D}, \mathbbm{L}, \mathbbm{U}$ stand for the diagonal, sub-diagonal and upper-diagonal parts of the matrix $\mathbbm{A}$, respectively. If the strict inequality $w_I < \tldvecu_I < v_I$
holds for some index $I$, then
%
%\begin{equation} zruseno cislovani - 3.5. Tomas
$$
\tldvecu_I=[ \left(1 - \omega \right)\tldvecu + \omega \mathbbm{D}^{-1}\left( 
\bfvec b - (\mathbbm{L}+\mathbbm{U})\tldvecu  \right)]_I.
$$
%\label{eq:minmax}
%\end{equation}
It means that 
\[
(\mathbbm{A} \tldvecu)_I = \bfvec b_I. 
\]
On the other hand, if $\tldvecu_I= w_I$ then $\tldvecu_I< v_I$. According to (\ref{eqA}) we conclude that 
\[
\tldvecu_I \ge [ \left(1 - \omega \right)\tldvecu + \omega \mathbbm{D}^{-1}\left( \bfvec b - (\mathbbm{L}+\mathbbm{U})\tldvecu  \right)]_I.
\]
Since $a_{II}>0$ we obtain 
\[
(\mathbbm{A} \tldvecu)_I \ge \bfvec b_I.
\]
Analogously, if $\tldvecu_I= v_I$ for some index $I$ then 
it follows from (\ref{eqA}) that 
\[
\tldvecu_I \le [ \left(1 - \omega \right)\tldvecu + \omega \mathbbm{D}^{-1}\left( \bfvec b - (\mathbbm{L}+\mathbbm{U})\tldvecu  \right)]_I.
\]
Thus
\[
(\mathbbm{A} \tldvecu)_I \le \bfvec b_I.
\]
Hence the vector $\tldvecu$ is a solution to the linear complementarity problem 
(\ref{eq:linear-problem-inequality}).

\section{Computational results}

\subsection{Convergence analysis in the case of smooth barrier functions}
\label{smooth}

In this section we present result of convergence analysis for a test problem in which the barrier functions are sufficiently smooth so that they belong to the space $W^{1,2}(\Omega)$.
We set up the following test problem. The computational domain is $[-0.5,0.5]^2$ and the initial condition is as follows:
\begin{equation}
u_{ini}\left(x,y\right)=\min\left\{x^2+y^2-0.55,0\right\}.
\label{eqn:eoc-initial-condition-ls}
\end{equation}
 The exact solution of (\ref{eq:level-set-segmentation})--(\ref{eq:level-set-segmentation-3}) with $g^0\equiv 1$ is 
 \begin{equation}
u\left(x,y,t\right)=\min\left\{x^2+y^2-0.55+t,0\right\}.
\label{eqn:exact-solution}
\end{equation}
%
%Experimental order of convergence of (\ref{eq:level-set-segmentation})--(\ref{eq:level-set-segmentation-3}) with very similar set-up has been studied in \cite{HandlovicovaMikulaOberhuber-2012}. 

\begin{figure}
\center{
\includegraphics[width=3cm,angle=-90]{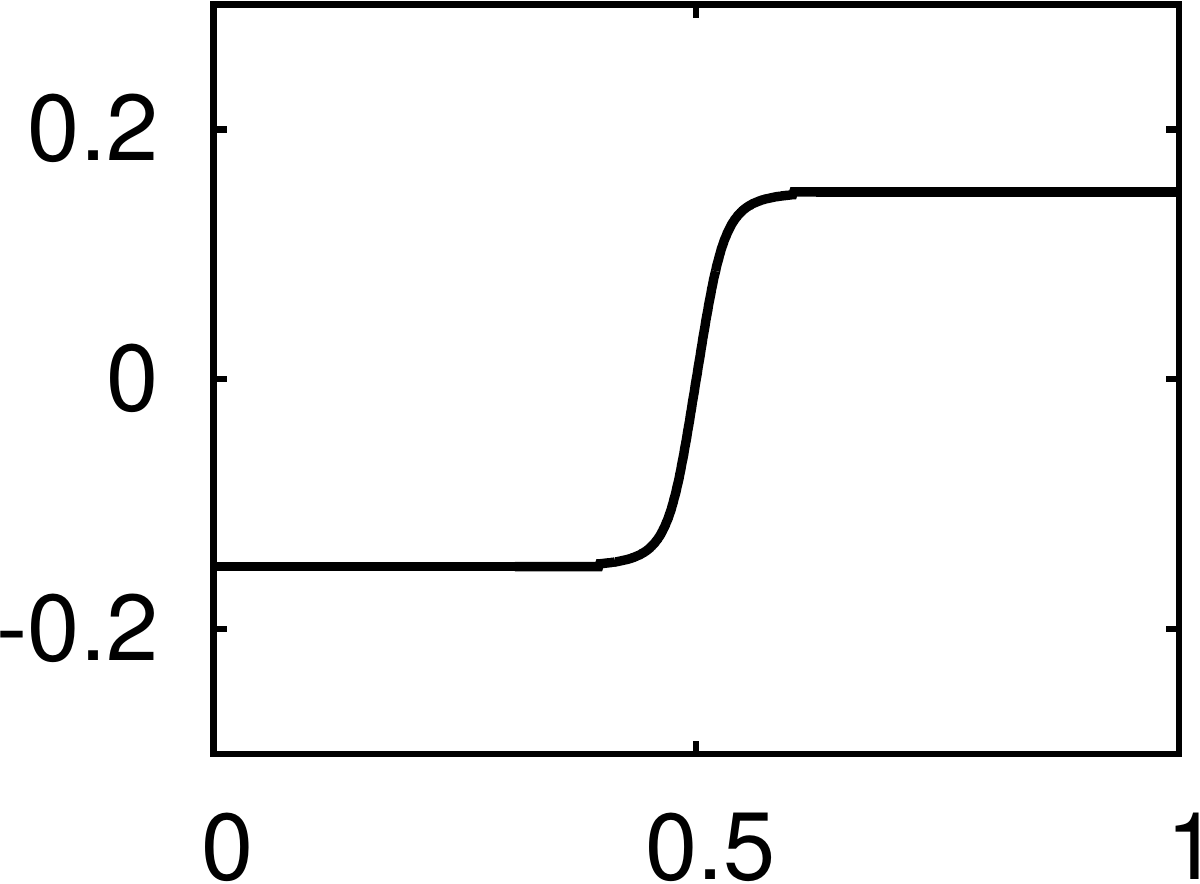}
\qquad 
\includegraphics[width=3cm,angle=-90]{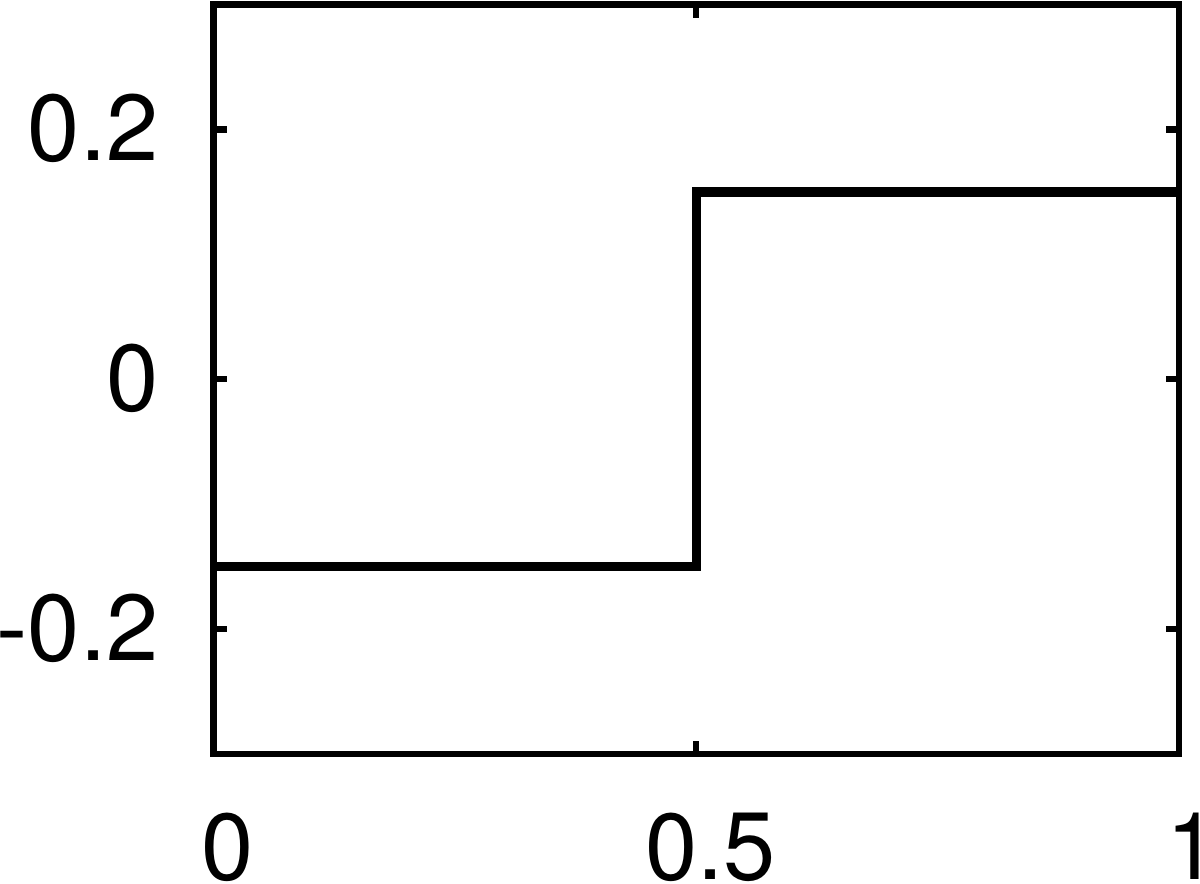}
%
% ze souboru data/bound.txt
%

\vspace{0.25cm}

\hspace{1cm} a) \hspace{5cm} b)
}
\caption{Different constraints considered for the convergence study: a) arctan shape, %b) rounded constraint with $r_0=0.15$, c) rounded constraint with $r_0=0.1$,
 b) discontinuous constraint.}
\label{fig:constraints}
\end{figure}

With additional smooth constraint $v \in W^{1,2}(\Omega)$, there is no longer an analytical solution. However, we see that the solution (\ref{eqn:exact-solution}) is axially symmetric with its center in origin. Prescribing an axially symmetric constraint enables us to keep axial symmetry of the solution. It allows us to transform our problem (\ref{eq:constrained-level-set-segmentation})--(\ref{eq:constrained-level-set-segmentation-3}) into coordinates $(r,t)$ where $r=\sqrt{x^2+y^2}$. The symmetric formulation then reads as
\begin{eqnarray}
\label{eq:symmetric-1}
&f_t = \frac{\epsilon^2}{\epsilon^2+\left(f_r^2 \right)} f_{rr} + \frac{f_r}{r}&\ \hbox{for}\ (r,t)\in(0,R)\times(0,T],\ \hbox{such that}\  w\left(r\right) < f\left(r\right) < v\left(r\right), \\
&f_t \geq \frac{\epsilon^2}{\epsilon^2+\left(f_r^2 \right)} f_{rr} + \frac{f_r}{r}&\ \hbox{for}\ (r,t)\in(0,R)\times(0,T],\ \hbox{such that}\  f\left(r\right) = w\left(r\right), \\
&f_t \leq \frac{\epsilon^2}{\epsilon^2+\left(f_r^2 \right)} f_{rr} + \frac{f_r}{r}&\ \hbox{for}\ (r,t)\in(0,R)\times(0,T],\ \hbox{such that}\  f\left(r\right) = v\left(r\right), \\
&f\left(\cdot,0\right) = f_{ini}&\ \hbox{for}\ r \in [0,R],\\
&f_r\left(R,t\right)=0&\ {\rm for}\ t \in [0,T], \\
\label{eq:symmetric-2}
&f_r\left(0,t\right)=0&\ \hbox{for}\ t \in [0,T].
\end{eqnarray}
In our numerical experiments, we set $w \equiv -\infty$. We solve the problem (\ref{eq:symmetric-1})--(\ref{eq:symmetric-2}) numerically with very fine resolution having $1000$ nodes on interval $[0,0.5]$. We take it as an exact solution for (\ref{eq:constrained-level-set-segmentation})--(\ref{eq:constrained-level-set-segmentation-3}) with symmetric initial condition $f_0\left(r\right)=\min\left\{r^2-0.55,0\right\}$ which agrees exactly with (\ref{eqn:eoc-initial-condition-ls}). Firstly, we set $v$ as 
\begin{equation}
v\left(r\right)= \left\{
   \begin{array}{cc}
   L/2 & {\rm for}\ r < R - r_0, \\
   \frac{L}{\pi} \arctan\left(\frac{r-R}{r_0}\left(5+35\left(\frac{\left(r-R\right)^2}{R^2}\right)\right)\right)  & {\rm for}\ R-r_0 \leq r \leq R + r_0, \\
   -L/2 & {\rm for}\ r > R + r_0. 
   \end{array}
   \right.
\label{eqn:arctan-like-constraint}   
\end{equation}
For $R=0.5$, $L=0.3$ and $r_0=0.1$ the shape of the function $v(r)$ is depicted on the Figure \ref{fig:constraints}.a.

\begin{figure}
\center{
\includegraphics[width=3cm,angle=-90]{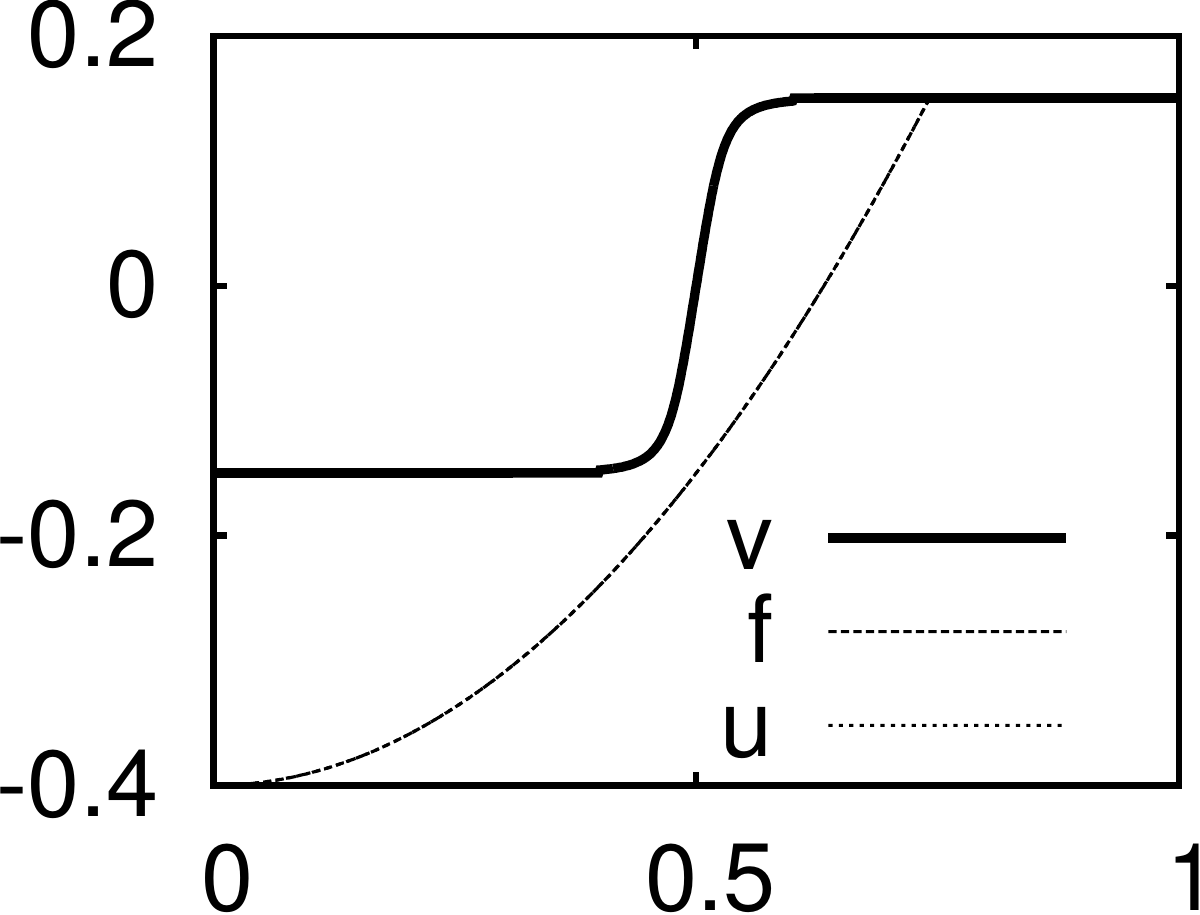} % ze souboru data/obstacle_arctg.txt data/sym_arctg_0.txt data/257_arctg_0.txt
\qquad
\includegraphics[width=3cm,angle=-90]{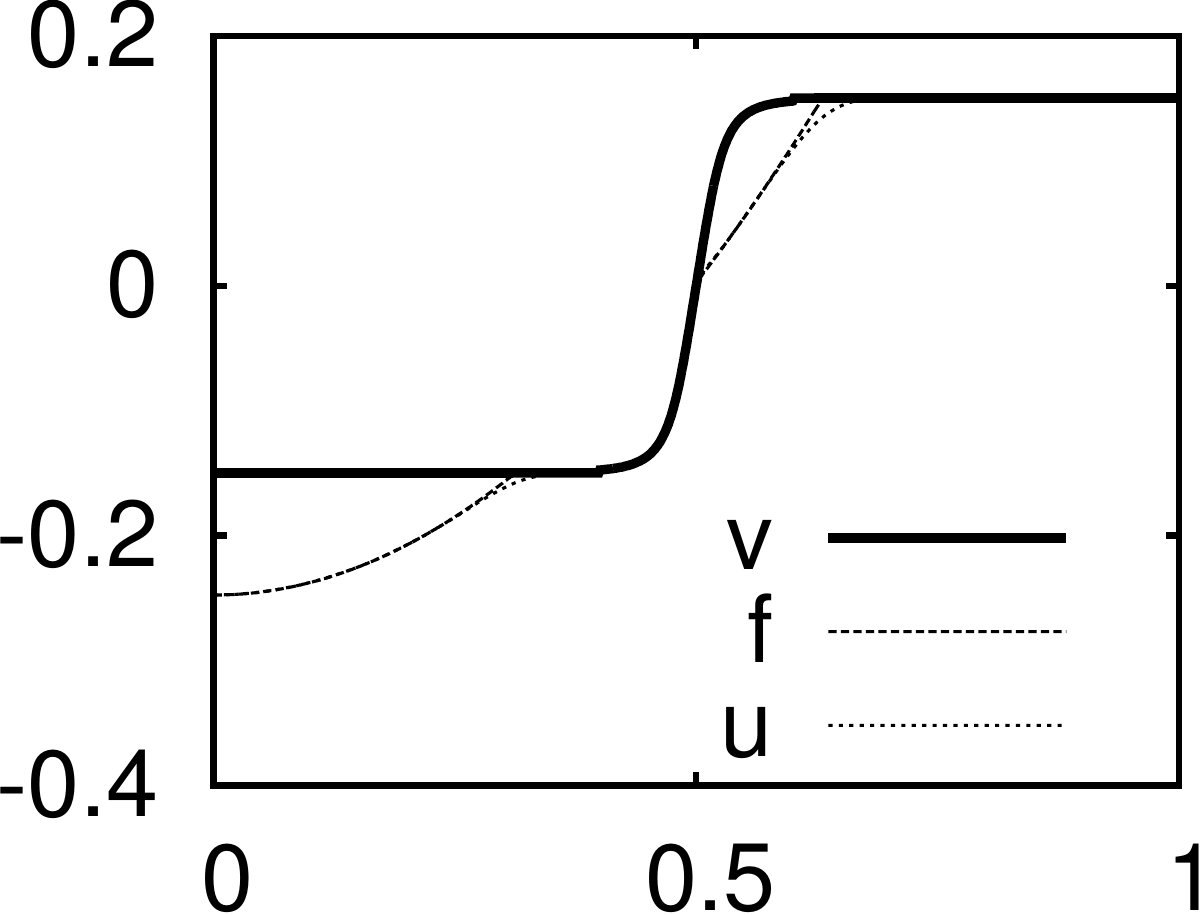} % se souboru data/obstacle_arctg.txt data/sym_arctg_0_075.txt data/257_arctg_0_075.txt

\vspace{0.25cm}

{\small \hspace{1cm} $t=0$ \hspace{4cm} $t=0.075$}

%\includegraphics[width=4cm,angle=-90]{fig3-3.pdf}

%\hspace{0.4cm} $t=0.15$ \hspace{3.5cm} $t=0.15$ (detail)
}
\caption{Results obtained with the arctan-like constraint $v$, $f$ is solution of symmetric problem (\ref{eq:symmetric-1}--\ref{eq:symmetric-2}) and $u$ is solution of level-set formulation (\ref{eq:constrained-level-set-segmentation})--(\ref{eq:constrained-level-set-segmentation-3}) with $g \equiv 1$.}
\label{fig:arctan-w}
\end{figure}

The results at times $t=0$ and $t=0.075$ are depicted on the Figure \ref{fig:arctan-w} and convergence rate analysis is summarized in Table \ref{tab:eoc-arctan-w}. We present the experimental order of convergence (EOC) for various norms of the error depending on the mesh discretization $h_1=h_2=1/N$ where $N=N_1=N_2$. It is defined as follows:
\[
 EOC = \log_2\left(\frac{error(h)}{error(h/2)}\right)
 \]

\begin{table}
\begin{center}
\begin{tabular}{l|l|l|l|l|l|l}\hline
\raisebox{-1ex}[0ex]{\footnotesize{$N=N_1=N_2$}}&
\multicolumn{2}{|c|}{\raisebox{1ex}[3.5ex]{$\left\| \cdot \right\|_{L^1\left(\omega_h \times \left(0,T\right)\right)}^{h,\tau}$}}&
\multicolumn{2}{|c|}{\raisebox{1ex}[3.5ex]{$\left\| \cdot \right\|_{L^2\left(\omega_h \times \left(0,T\right)\right)}^{h,\tau}$}}&
\multicolumn{2}{|c}{\raisebox{1ex}[3.5ex]{$\left\| \cdot \right\|_{L^\infty\left(\omega_h \times \left(0,T\right)\right)}^{h,\tau}$}} \\
 \cline{2-7}
         &  \footnotesize{Error} &      {\footnotesize EOC}        & \footnotesize{Error}        & \footnotesize{EOC}& \footnotesize{Error}&      \footnotesize{EOC}\\ \hline \hline 
\footnotesize{$32$}   & 0.004232  &                            & 0.000752    &                           & 0.087758  &                        \\
\footnotesize{$64$}   & 0.001777  &    \raisebox{1ex}{\bf 1.25}& 0.000174    & \raisebox{1ex}{\bf 2.11}  & 0.046046  &    \raisebox{1ex}{\bf 0.93}\\
\footnotesize{$128$}  & 0.000732  &    \raisebox{1ex}{\bf 1.27}& 0.000044    & \raisebox{1ex}{\bf 1.98}  & 0.0249    &    \raisebox{1ex}{\bf 0.88}\\
\footnotesize{$256$}  & 0.000328  &    \raisebox{1ex}{\bf 1.15}& 0.000012    & \raisebox{1ex}{\bf 1.87}  & 0.015742  &    \raisebox{1ex}{\bf 0.66}\\
\footnotesize{$512$}  & 0.000152  &    \raisebox{1ex}{\bf 1.10}& 0.000003    & \raisebox{1ex}{\bf 2   }  & 0.01068   &    \raisebox{1ex}{\bf 0.55}\\
\hline
\end{tabular}
\caption{Experimental order of convergence with the arctan-like constraint $v$ given by (\ref{eqn:arctan-like-constraint}). It shows the second order of convergence in the $L^2$ norm while the convergence rate in the $L^1$ is only linear in $1/N$. }
\label{tab:eoc-arctan-w}
\end{center}
\end{table}

\subsection{Decrease of convergence ratio for non-smooth barrier functions}
\label{nonsmooth}

\begin{figure}
\center{
\includegraphics[width=3cm,angle=-90]{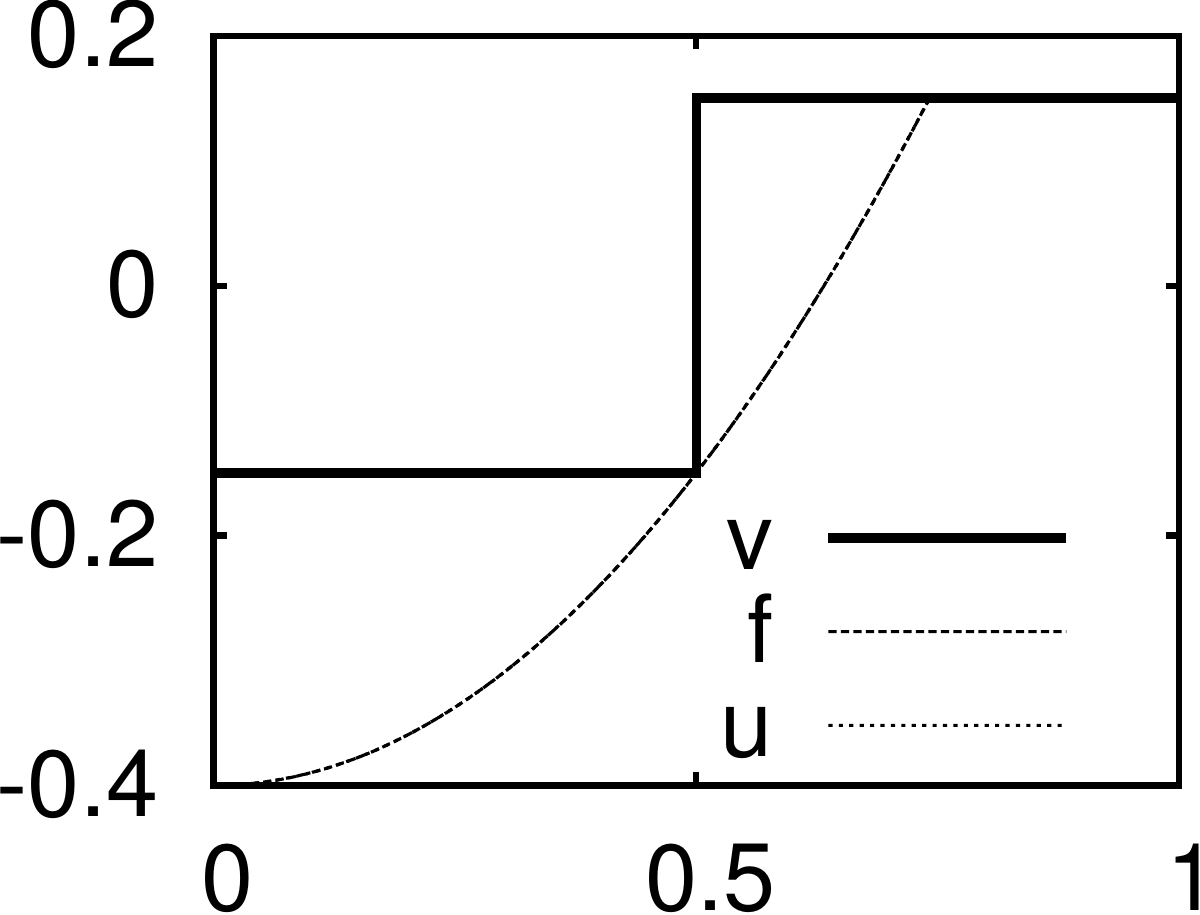} % ze souboru data/bound.txt data/symetric_t=0.txt data/computed_t=0.txt
\qquad
\includegraphics[width=3cm,angle=-90]{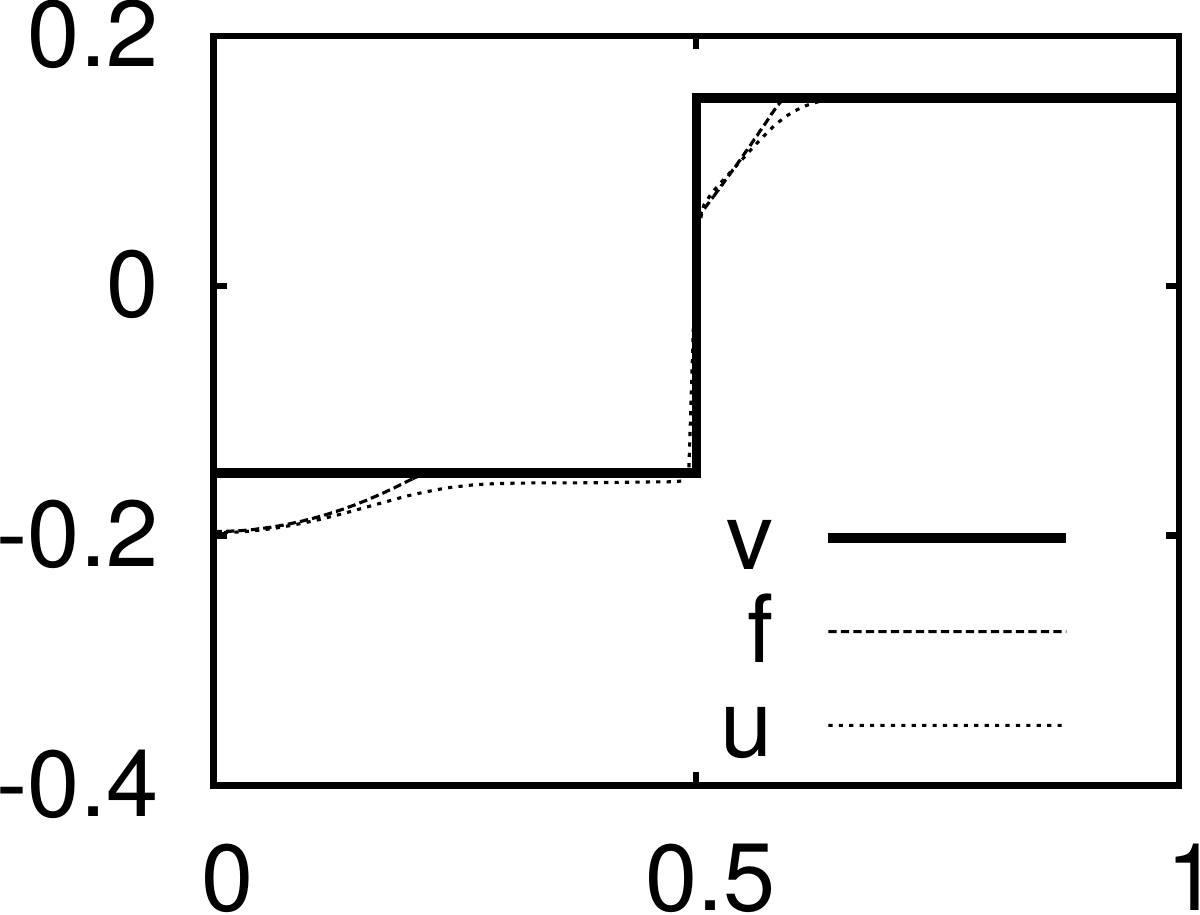} % ze souboru data/bound.txt data/symetric_t=0.101532.txt data/computed_t=0.101532.txt
\qquad 
\includegraphics[width=3cm,angle=-90]{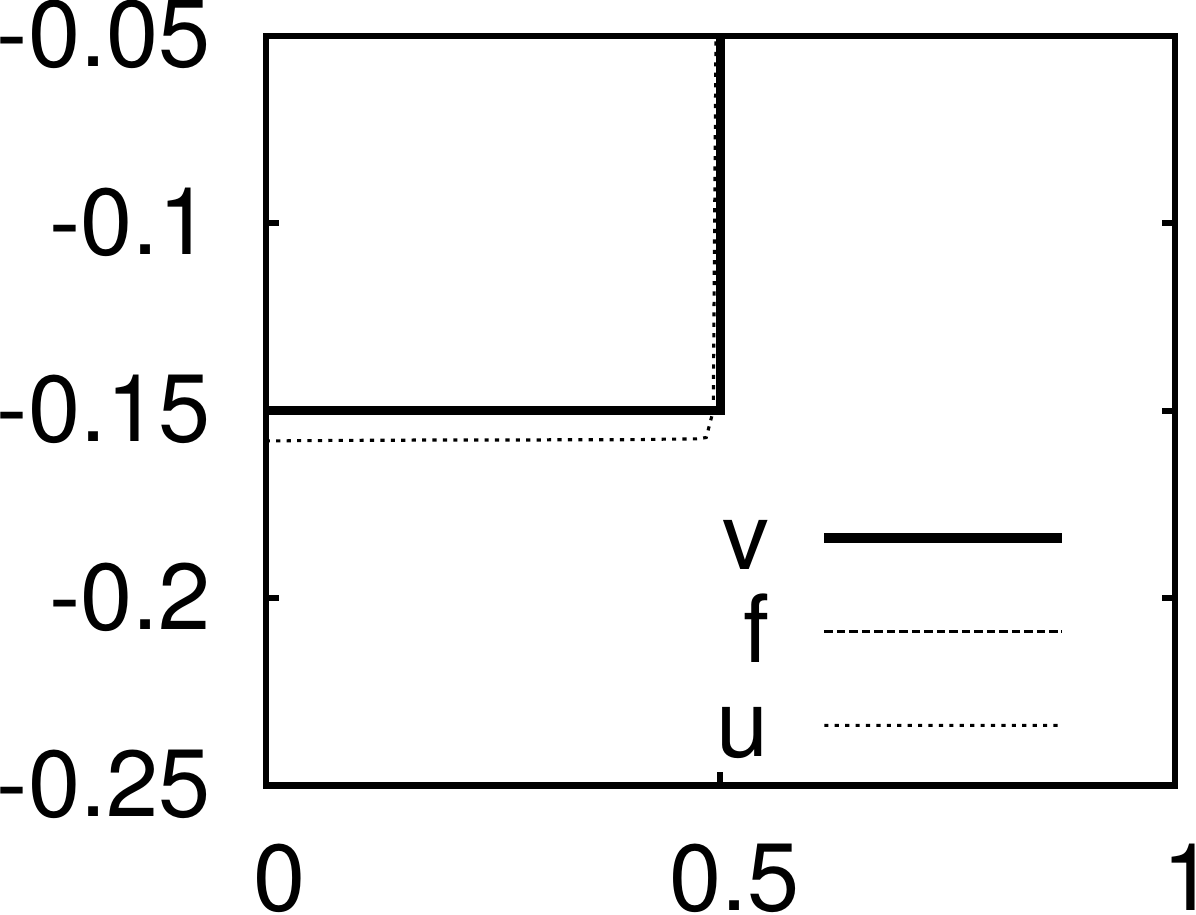} % ze souboru data/bound.txt data/symetric_t=0.15234.txt data/computed_t=0.15234.txt

\vspace{0.25cm}

{\small \hspace{2cm} $t=0$ \hspace{3cm} $t=0.1$ \hspace{3cm} $t=0.15$ (detail)}
}

\vspace{0.25cm}

\caption{Results obtained with the discontinuous constraint $v$, $f$ is solution of symmetric problem (\ref{eq:symmetric-1}--\ref{eq:symmetric-2}) and $u$ is solution of level-set formulation (\ref{eq:constrained-level-set-segmentation})--(\ref{eq:constrained-level-set-segmentation-3}) with $g \equiv 1$.}
\label{fig:discontinous-w}
\end{figure}

In this section, we present an example showing that the assumption of smoothness of barrier functions has a strong impact on the convergence of the method. We consider an upper barrier function $v$ representing a discontinuous constraint:
\begin{equation}
v\left(r\right) = \left\{
   \begin{array}{cc}
   -0.15 & {\rm for}\ r \leq 0.5, \\
    0.15   & {\rm for}\ r \geq 0.5.
   \end{array}
   \right.
\end{equation}
This is the most common choice in practical image segmentation computations. However, such a function does not belong to the space  $W^{1,2}(\Omega)$. The shape of the constraint is depicted on Figure~\ref{fig:constraints}.b and the results are summarized in Table \ref{tab:eoc-discontinous-w}. Details of the radially symmetric solution profile are depicted on the Figure~\ref{fig:discontinous-w}. We see that a thin interface develops between the level-set function $u$ and the constraint $v$. It is reflected by the $L^\infty$ norm of the error in which the scheme converges very slowly. The reason why the convergence is slowed is due to the fact that the set of discontinuity of a barrier function is fixed. As a consequence, a solution to a variational inequality is pasted to the barrier function in this set of discontinuity.

\begin{table}
\begin{center}
\begin{tabular}{l|l|l|l|l}\hline
\raisebox{-1ex}[0ex]{\footnotesize{$N=N_1=N_2$}}&
\multicolumn{2}{|c|}{\raisebox{1ex}[3.5ex]{$\left\| \cdot \right\|_{L^1\left(\omega_h \times \left(0,T\right)\right)}^{h,\tau}$}}&
\multicolumn{2}{|c}{\raisebox{1ex}[3.5ex]{$\left\| \cdot \right\|_{L^2\left(\omega_h \times \left(0,T\right)\right)}^{h,\tau}$}} \\
%\multicolumn{2}{|c|}{\raisebox{1ex}[3.5ex]{$\left\| \cdot \right\|_{L^\infty\left(\omega_h;\left[0,T\right]\right)}^{h,\tau}$}} \\
 \cline{2-5}
         &  \footnotesize{Error} &      {\footnotesize EOC}        & \footnotesize{Error}        & \footnotesize{EOC} \\ \hline \hline % \footnotesize{Error}&      \footnotesize{EOC}\\ \hline \hline 
\footnotesize{$32$}   & 0.003713  &                            &   0.000777  &                          \\ %&  0.142039 &                           \\
\footnotesize{$64$}   & 0.002179  &    \raisebox{1ex}{\bf 0.77}&   0.000256  & \raisebox{1ex}{\bf 1.6}  \\ %&  0.161288 &    \raisebox{1ex}{\bf -0.18}\\
\footnotesize{$128$}  & 0.001287  &    \raisebox{1ex}{\bf 0.76}&   0.000096  & \raisebox{1ex}{\bf 1.4}  \\ %&  0.178583 &    \raisebox{1ex}{\bf -0.14}\\
\footnotesize{$256$}  & 0.000823  &    \raisebox{1ex}{\bf 0.64}&   0.00004   & \raisebox{1ex}{\bf 1.3}  \\ %&  0.18882  &    \raisebox{1ex}{\bf -0.08}\\
\hline
\end{tabular}
\caption{Experimental order of convergence with the discontinuous constraint $v$.}
\label{tab:eoc-discontinous-w}
\end{center}
\end{table}

\subsection{Smoothing of the barrier functions}
\label{smoothing}
In previous parts \ref{smooth} and \ref{nonsmooth} we showed the importance of a smooth barrier function and its impact on the accuracy of computed solutions. If the barrier is smooth than the accuracy of the level set solution is higher. However, in practice the barrier functions $w$ and $v$ are often characteristic functions of the prescribed regions in the images. Nevertheless, we can smooth such step functions in a canonical way using their convolution with a two dimensional Gaussian kernel $G_\sigma(x)=(2\pi\sigma^2)^{-1} \exp(-\Vert x\Vert^2/2\sigma^2)$ where $\sigma<<1$ is sufficiently small. It means that the barrier functions are defined as:
\[
w_\sigma = G_\sigma * w, \qquad v_\sigma = G_\sigma * v.
\]
The advantage of a such smoothing is the order preserving property, i.e. $w(x)<v(x)$ for each $x\in\Omega$ implies $w_\sigma(x)<v_\sigma(x)$ in $\Omega$.

\section{Application to image segmentation and computational results}
\label{sect:results}
\begin{figure}
\center{
\includegraphics[width=4cm]{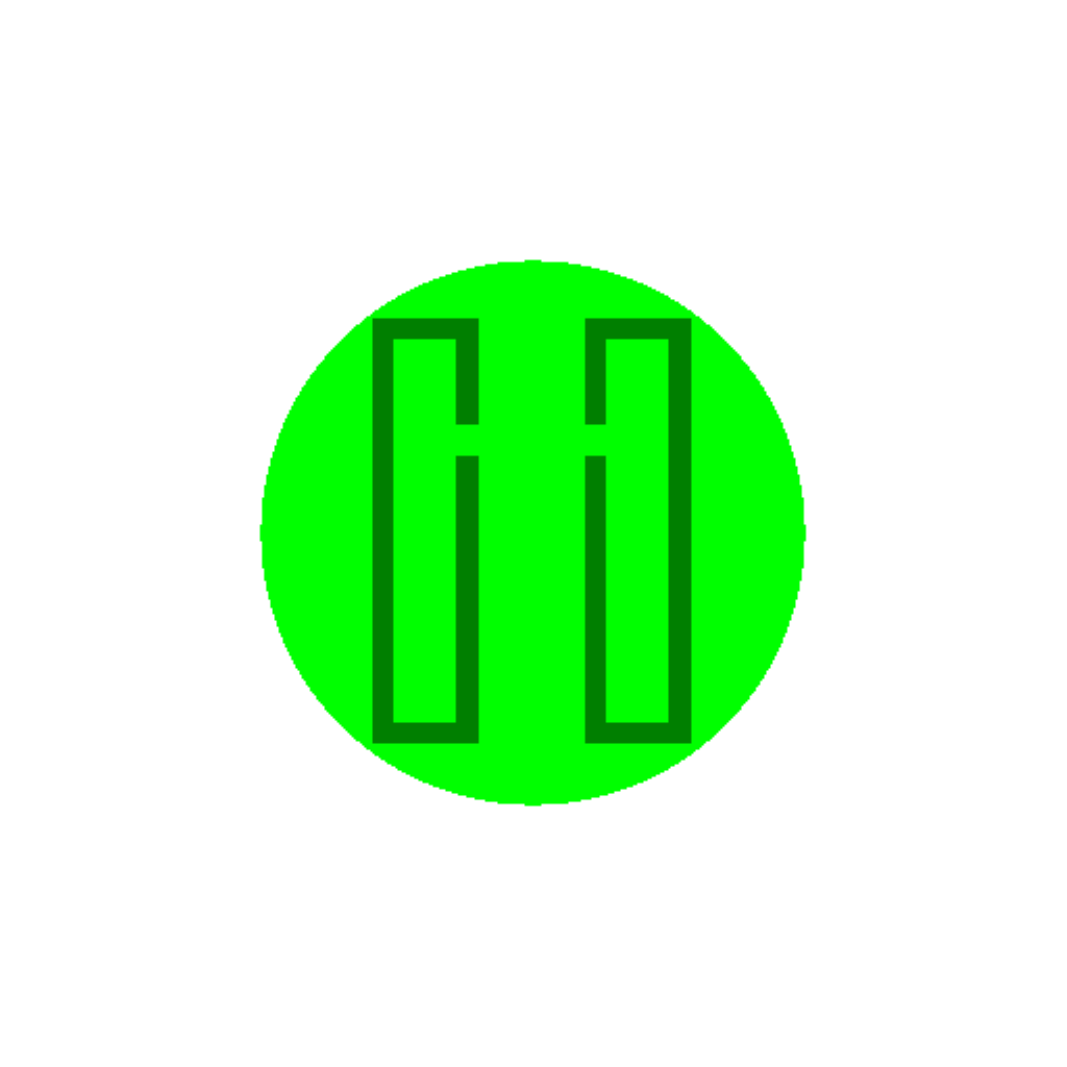}
\includegraphics[width=6cm]{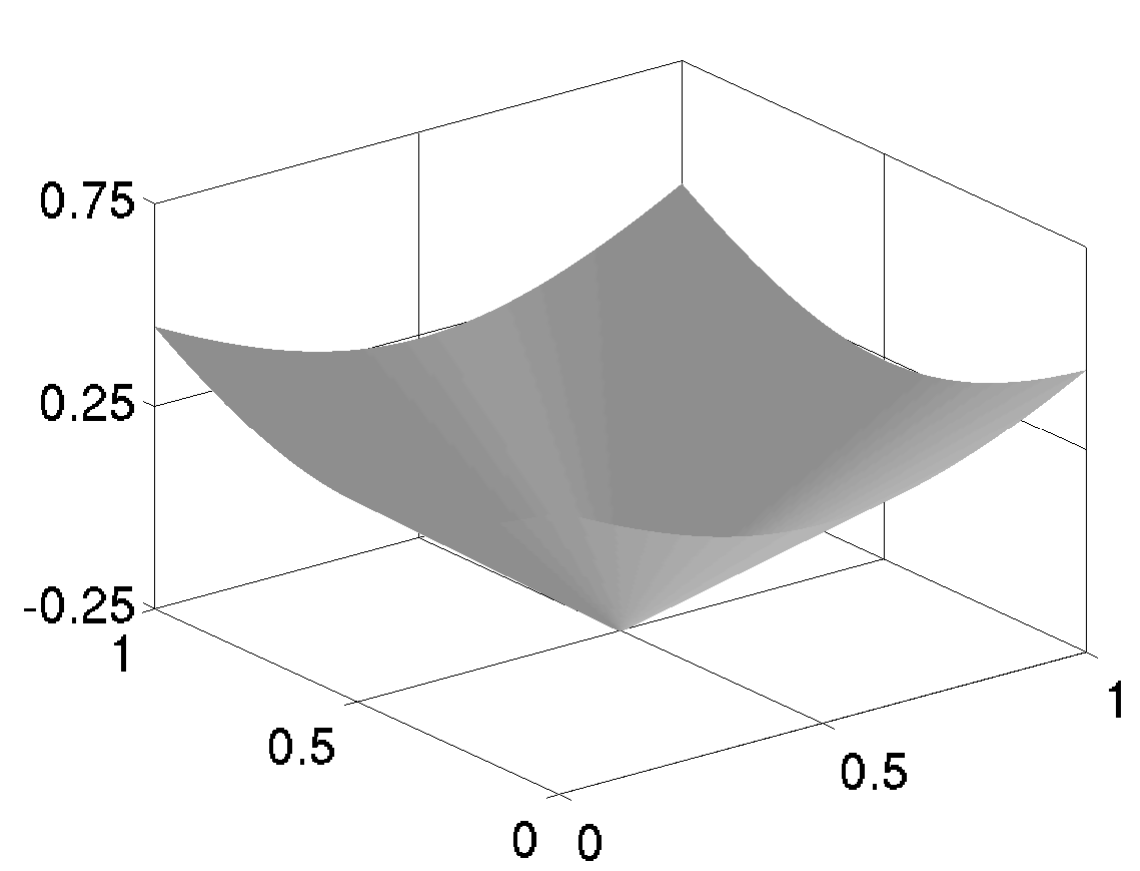}
\centerline{a)}
\includegraphics[width=4cm]{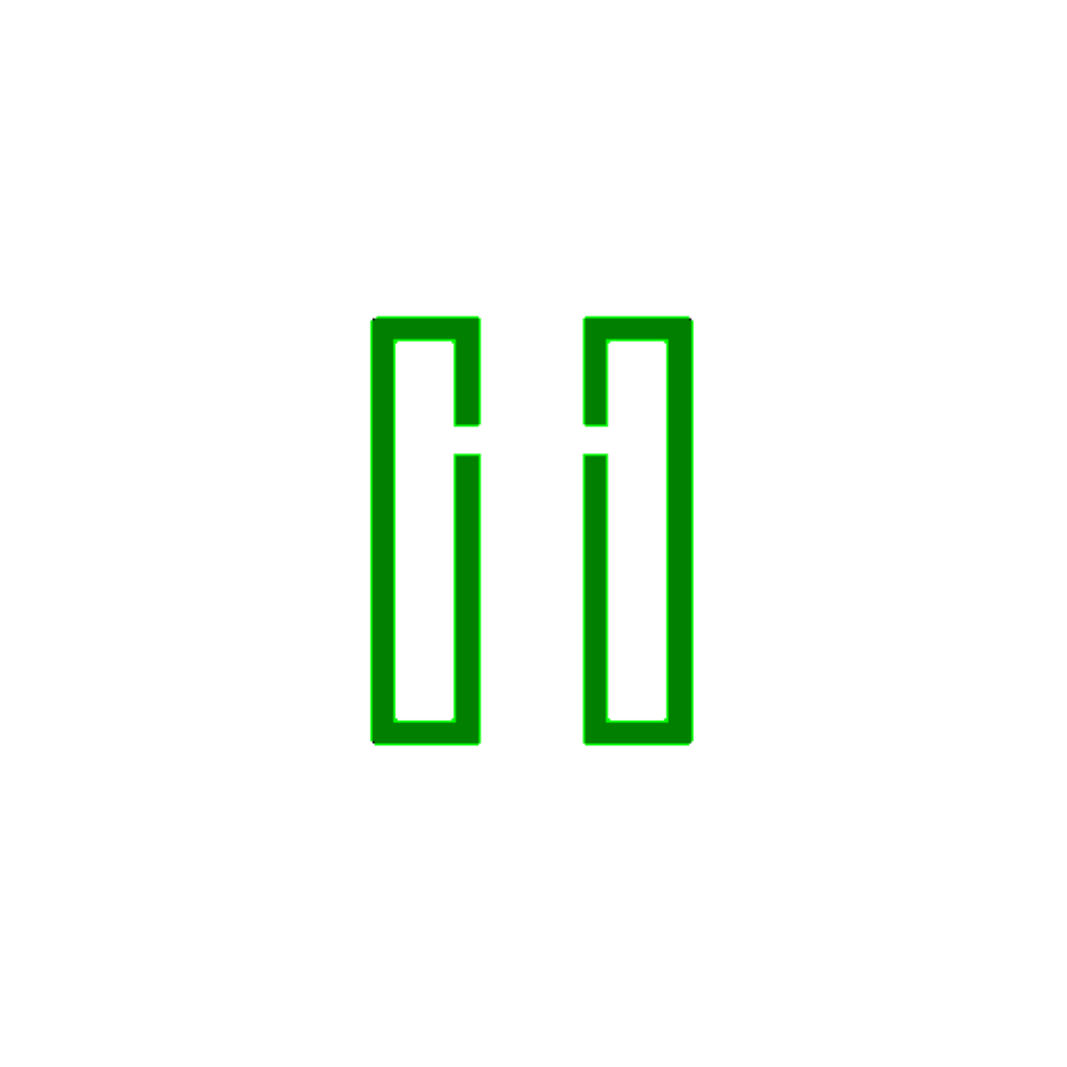} % cas = 1.2
\includegraphics[width=6cm]{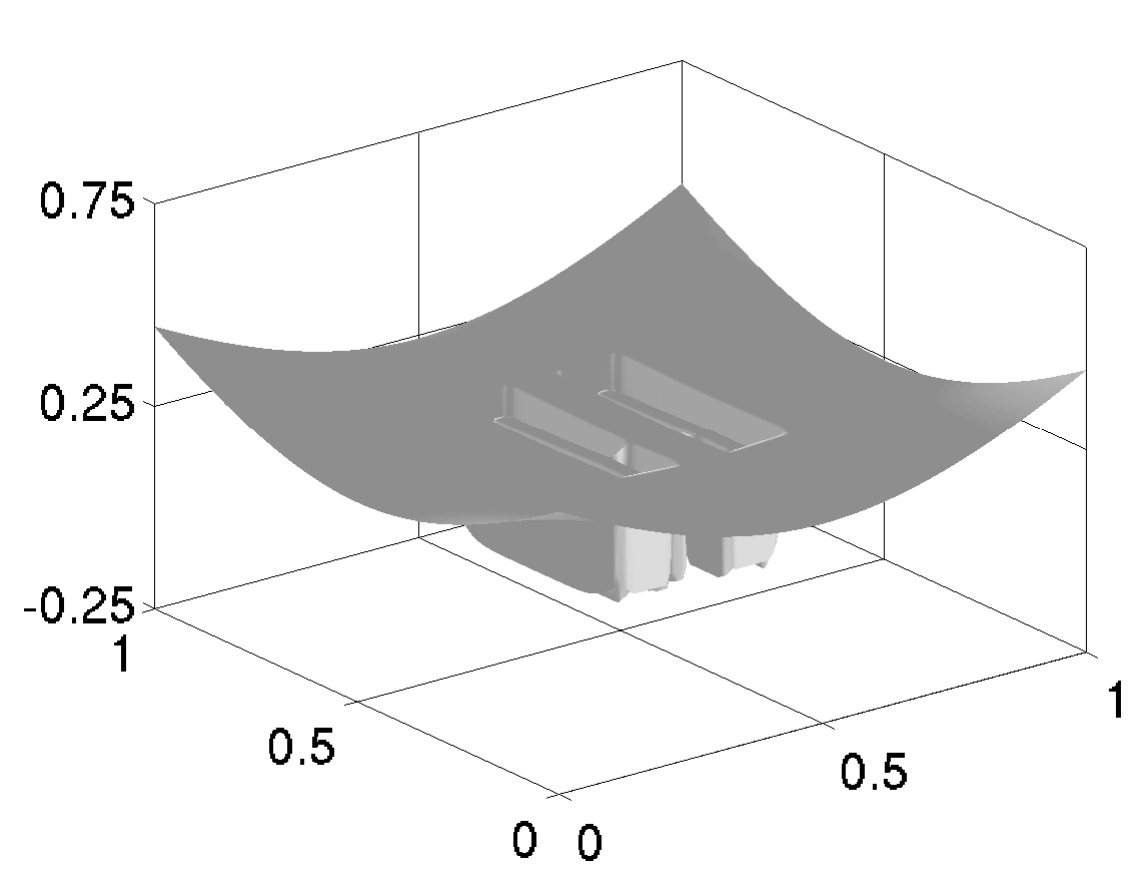} % cas = 1.2
\centerline{b)}
\includegraphics[width=4cm]{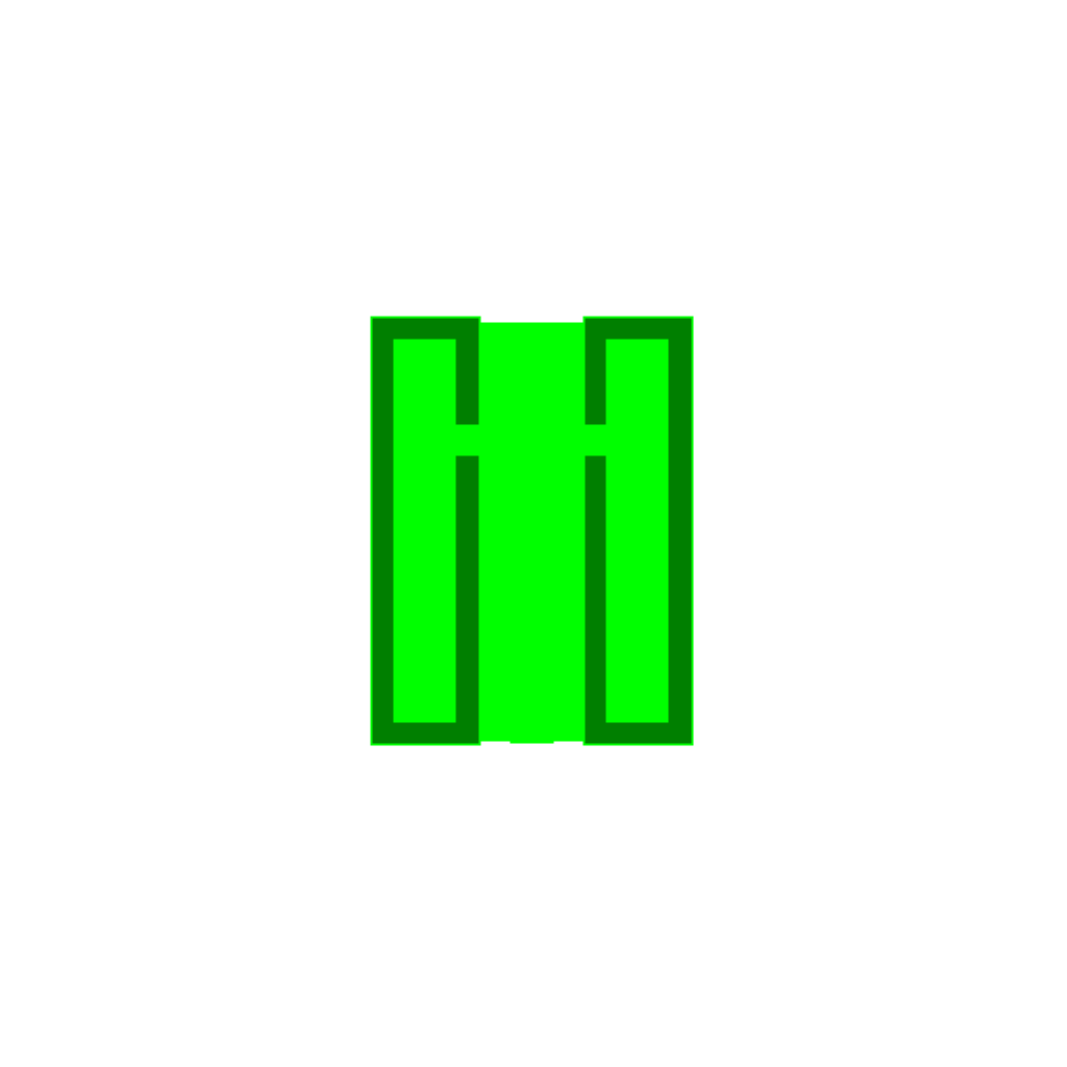} % cas = 145
\includegraphics[width=6cm]{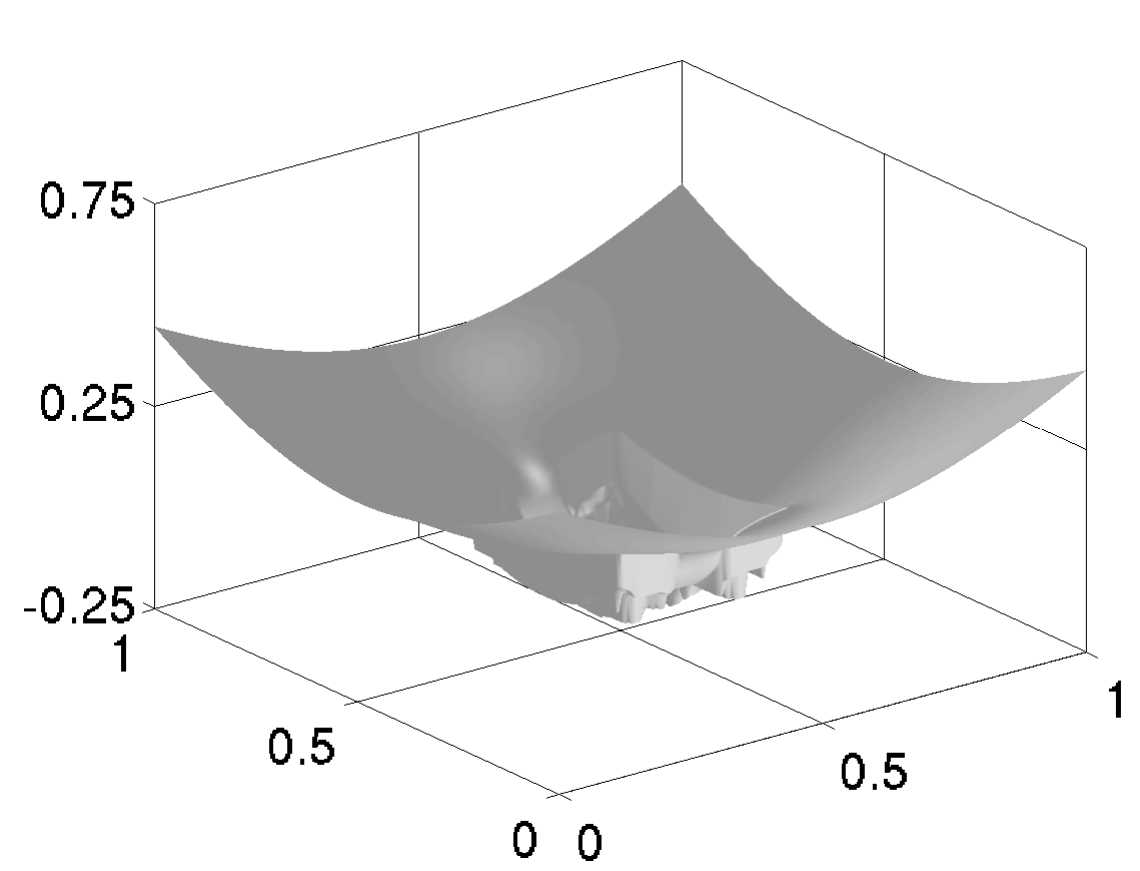} % cas = 145 eps = 0.0001
\centerline{c)}
}
\caption{a): We plot an initial curve $\Gammanul$ set as a circle with the radius $r=\sqrt{0.08}$ (left) and its level-set function (right). b): The segmentation at time $t=1.2$ after setting the regularizing parameter $\epsilon$ from (\ref{def:Q}) to $1$ is depicted together with the level-set function. c): The segmentation result with $\epsilon=0.0001$ at time $t=145$.}
\label{fig:no-constraint}
\end{figure}

In this section, we present experimental results obtained by the constrained level-set method. We first demonstrate the effect of the constraints applied to artificial images. Fig.~\ref{fig:no-constraint} a) (left) shows an image we want to segment. There are two rectangles with centers at points $(0.4,0.5)$ and $(0.6,0.5)$. The width of the rectangles is $0.1$ and the height is $0.4$. Thickness of the rectangles edges is $0.04$. On the inner edge of each rectangle, there is a thin hole. The image intensity function is defined on the domain $\Omega \equiv (0,1)^2$. The initial curve is a circle with the radius $r=\sqrt{0.08}$ centered exactly between the rectangles. Its signed distance function (taken as the initial level-set function $u_{ini}$) is depicted in Fig.~\ref{fig:no-constraint} a) (right). The numerical mesh consisted of $128\times 128$ grid points. If the thin holes in the rectangles should be taken into account we may want to segment only the edges of rectangles. It can be achieved by setting the regularizing parameter $\epsilon=1$. In \cite{MikulaSarti-2007} Mikula and Sarti interpreted the regularization parameter $\epsilon$ as a parameter by which we can control "convexity" of the final segmentation curve. The result taken at time $t=1.2$ is depicted in Fig.~\ref{fig:no-constraint} b). On the other hand, if small holes are present by mistake, we may want the segmentation curve to fill them. It may be achieved by setting $\epsilon=0.0001$. However, a convex hull of both rectangles is segmented as we can see in Fig.~\ref{fig:no-constraint} c).

\begin{figure}
\center{
\includegraphics[width=5cm]{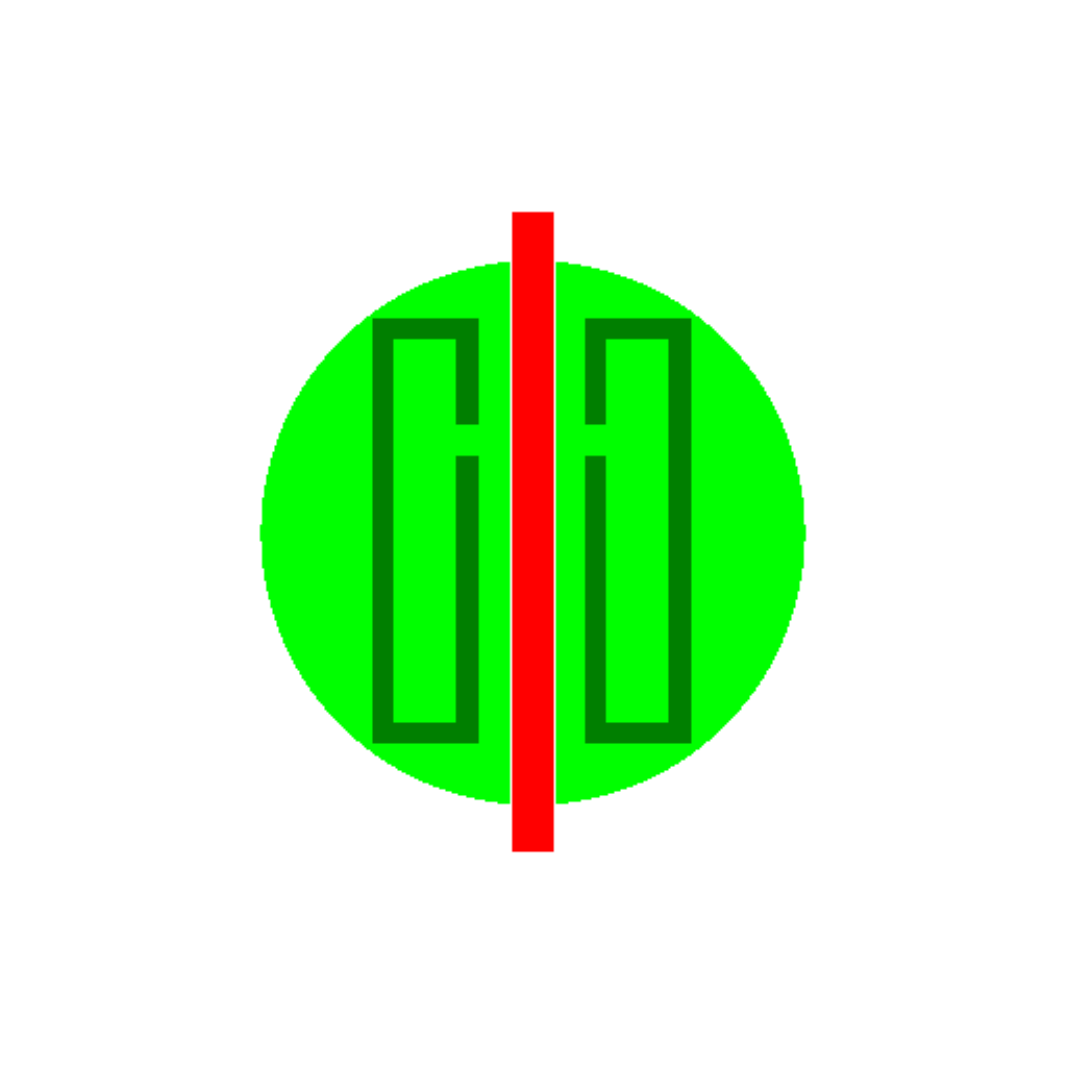}
\includegraphics[width=7cm]{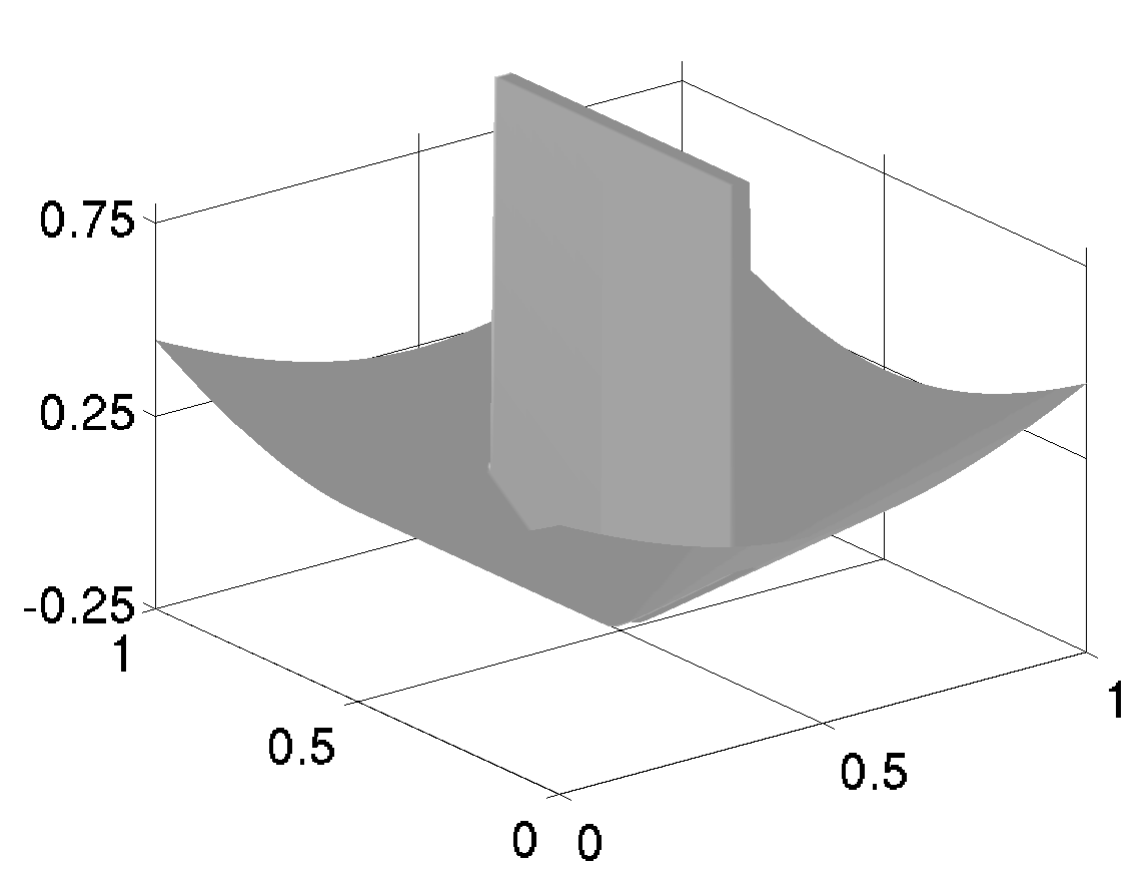}
\centerline{a)}
\includegraphics[width=5cm]{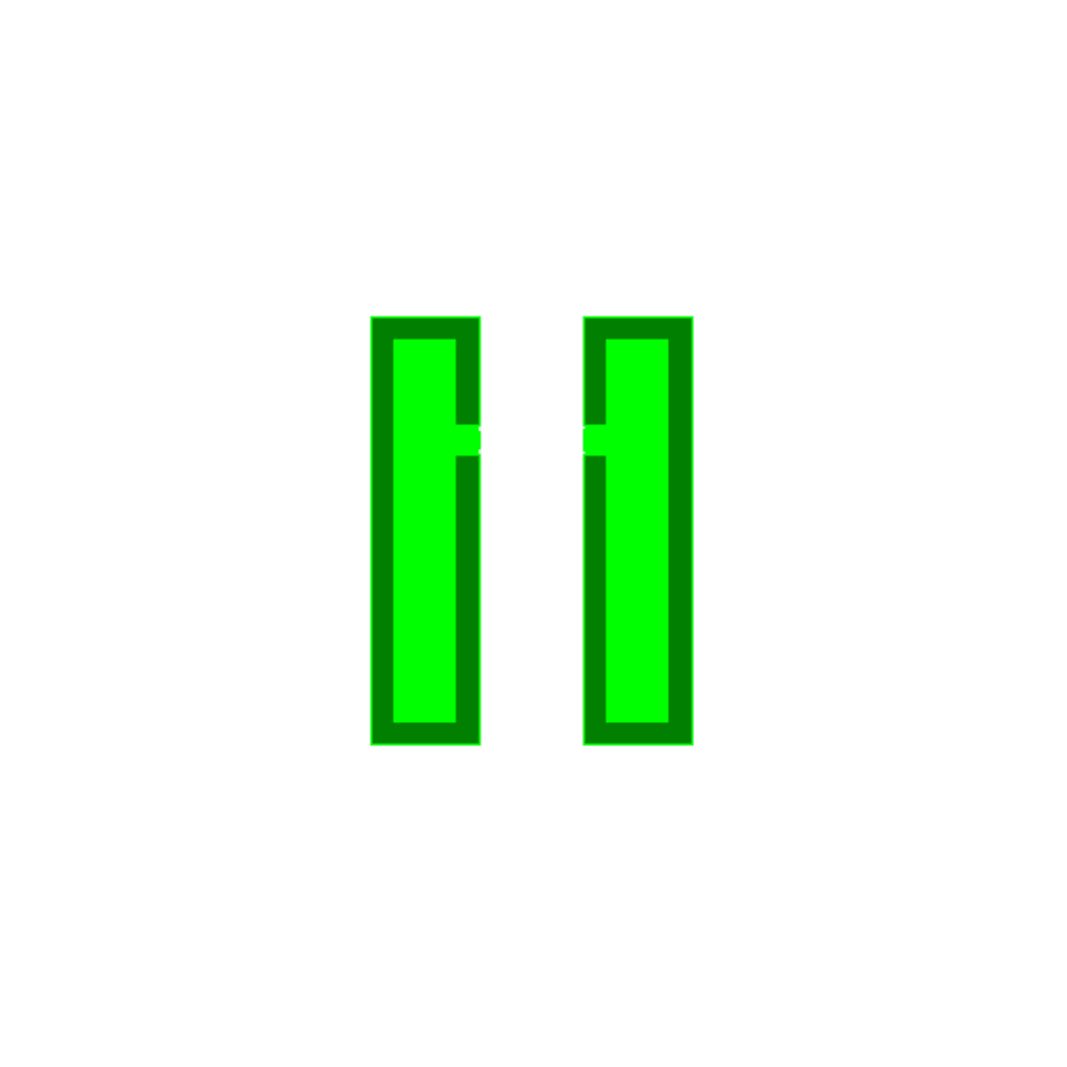} % cas = 67 eps = 0.0001
\includegraphics[width=7cm]{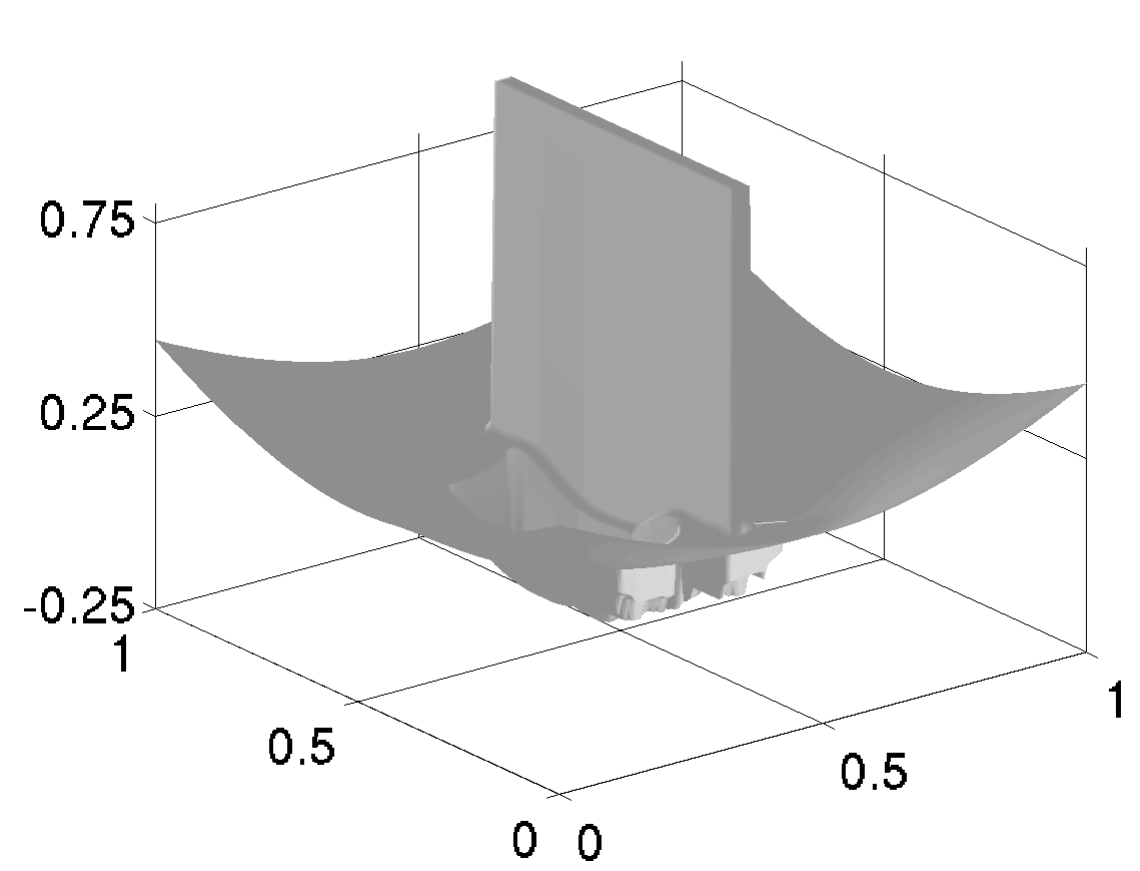} % cas = 67 eps = 0.0001
\centerline{b)}
}
\caption{a): An example of prescription of a single constraint depicted by the red bar in the middle of the image. The initial curve $\Gammanul$ is a circle with the radius $r=\sqrt{0.08}$. b): The segmentation at time $t=67$ shows nice separation of both rectangles. The corresponding level-set functions are depicted on the right.}
\label{fig:one-constraint}
\end{figure}

The segmentation we aim to can be achieved by prescribing one constraint guaranteeing that the part of the image between the rectangles must be outside the segmentation domain. We construct the constraint by putting a red bar placed on between the rectangles on the Fig.~\ref{fig:one-constraint} a). Its width equals $0.04$ and the height is set to $0.6$. The constraint function $w(x)$ is positive on the red bar region and negative everywhere else. We again set $\epsilon=0.0001$. The result taken at time $t=67$ is shown on Fig.~\ref{fig:one-constraint} b).

\begin{figure}
\center{
\includegraphics[width=5cm]{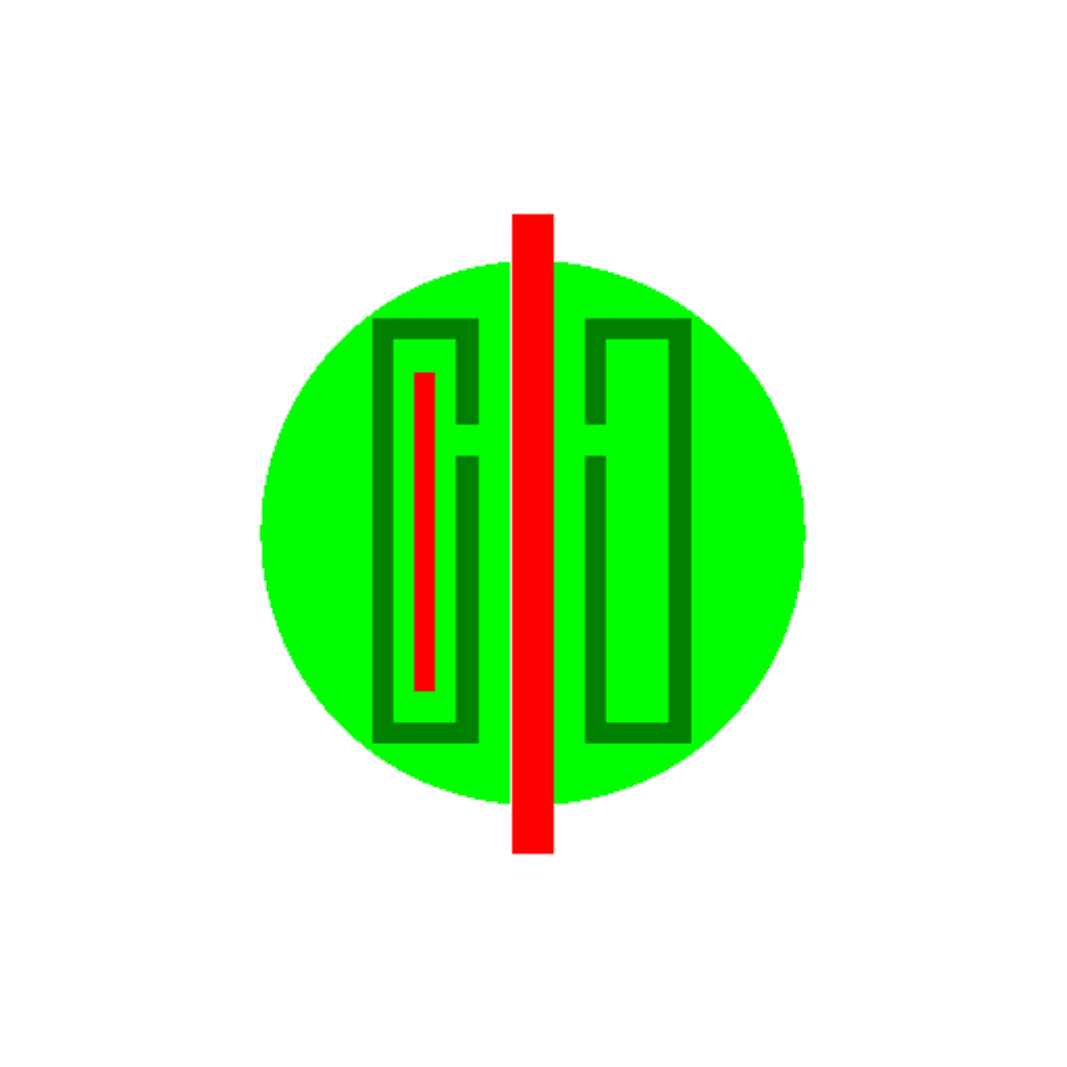}
\includegraphics[width=7cm]{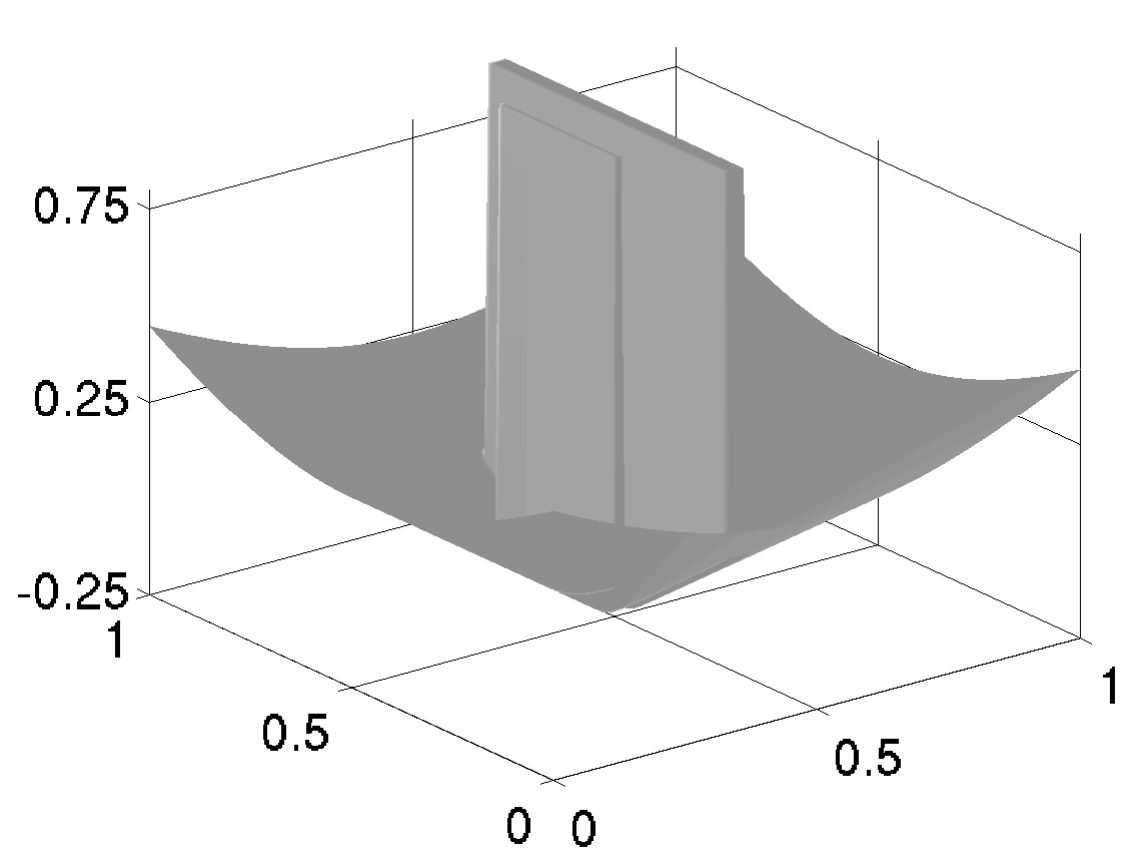}
\centerline{a)}
\includegraphics[width=5cm]{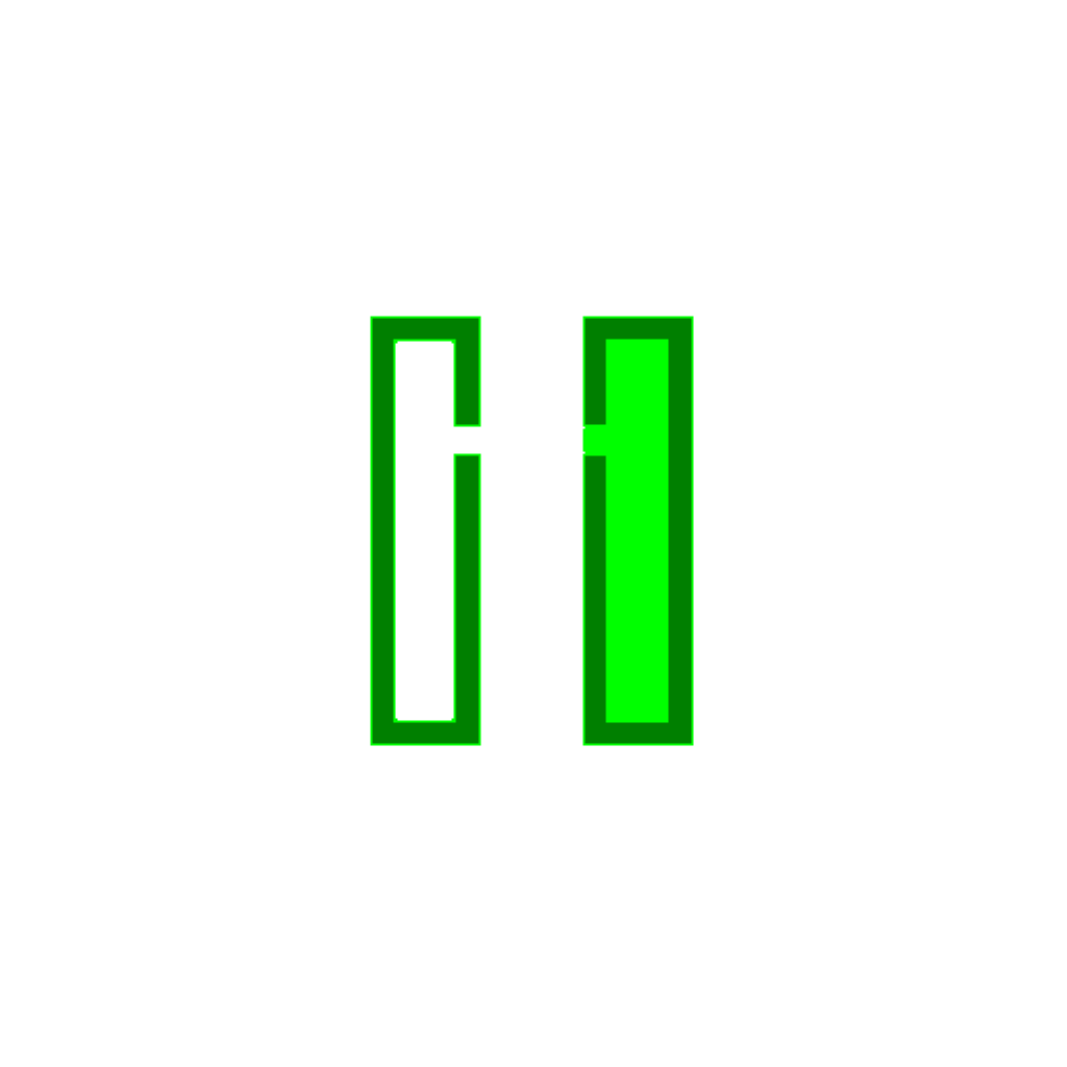} % cas = 67 eps = 0.0001
\includegraphics[width=7cm]{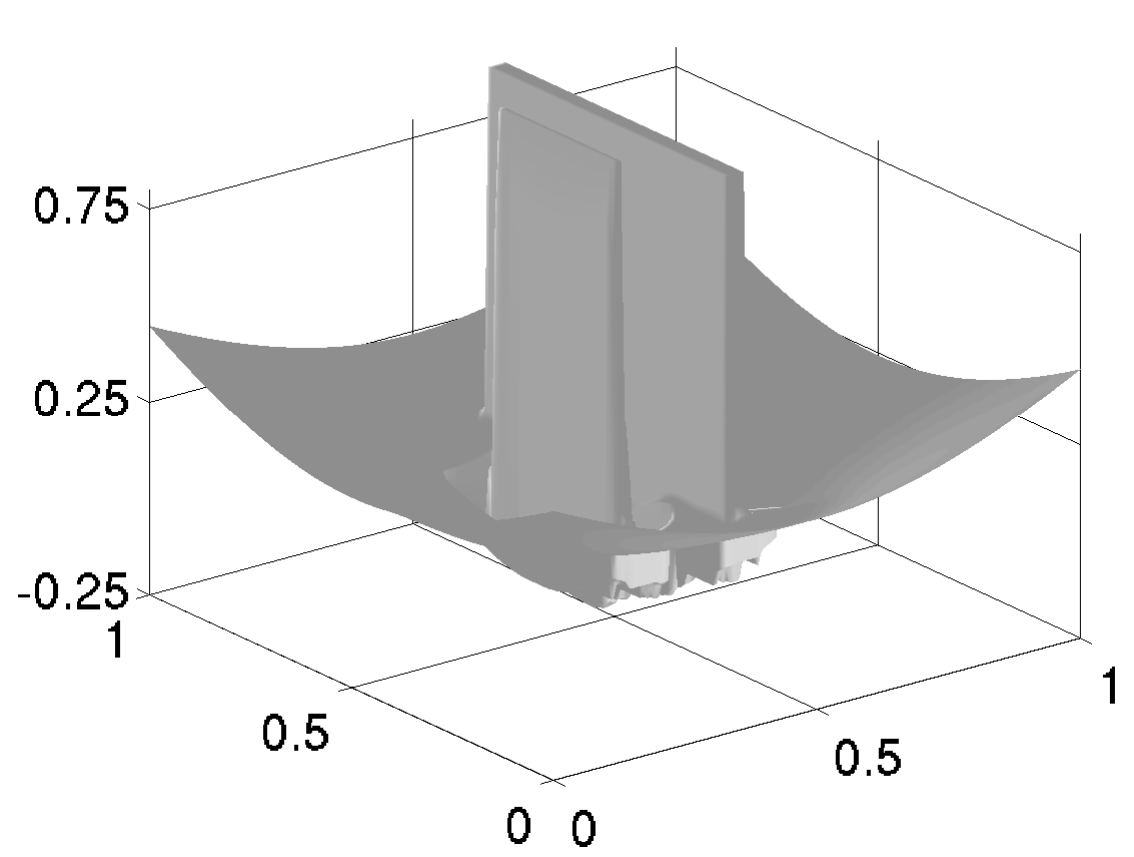} % cas = 67 eps = 0.0001
\centerline{b)}
}
\caption{a): An example with two prescribed constraints depicted by red bars (top left). The initial curve $\Gammanul$ is a circle with the radius $r=\sqrt{0.08}$. b): The segmentation result at time $t=67$ is shown in the bottom left part of the figure. The corresponding level-set functions are shown on the right.}
\label{fig:two-constraints}
\end{figure}

In our last synthetic example shown in Fig.~\ref{fig:two-constraints} we place one more constraint inside the left rectangle. The segmentation poses features of both experiments on Figures \ref{fig:no-constraint} and \ref{fig:one-constraint}. It shows that the constrained level-set method is capable of controlling which part of the image we want to segment.

\begin{figure}
\center{
\includegraphics[width=4cm]{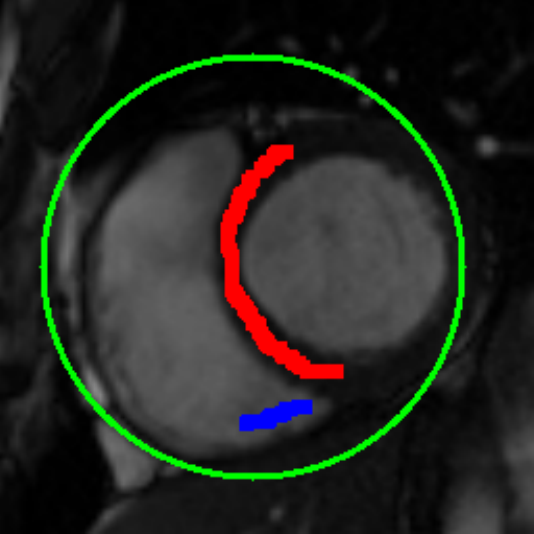}
\includegraphics[width=4cm]{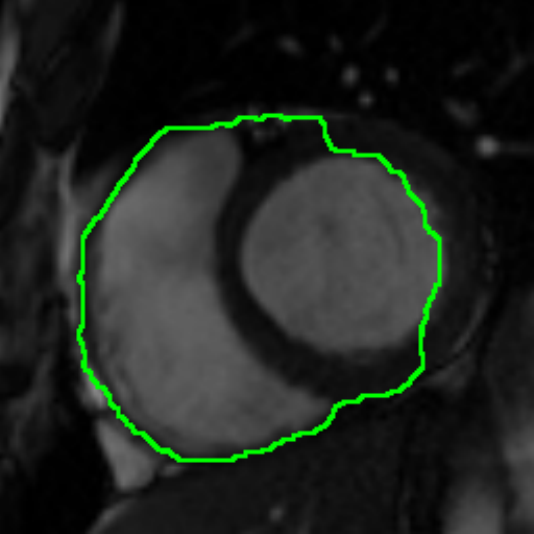}
\includegraphics[width=4cm]{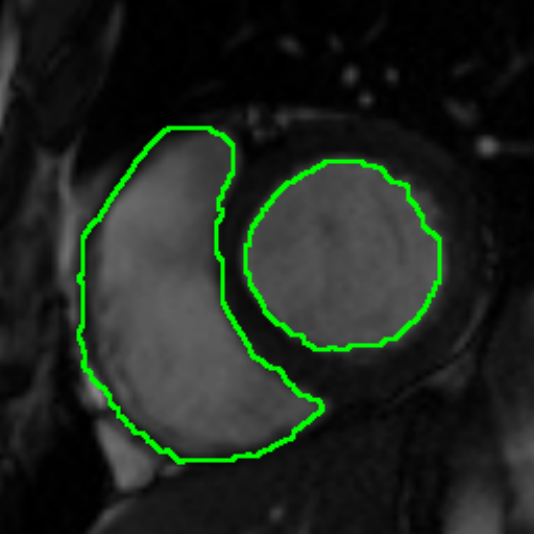}
\centerline{a) \hspace{4cm} b) \hspace{4cm} c)} \\
\includegraphics[width=4cm]{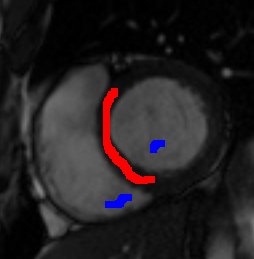}
\includegraphics[width=4cm]{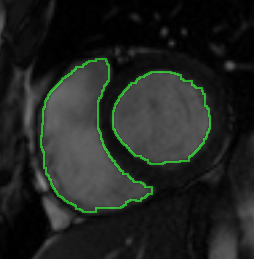}
\centerline{d) \hspace{4cm} e)} 
}
\caption{a) Setup for the constrained level set method with the initial curve $\Gammanul$ (green) and two constraints -- the red curve stands for the exterior while the blue one is for interior of the segmented region. b) Segmentation results obtained by the unconstrained level-set method. c) Segmentation obtained by means of the constrained level-set method. The interior constraint helps to capture the bottom of the right ventricle. d) Initial seeds for the graph-cuts method. e) Segmentation obtained by the graph-cuts method.}
\label{fig:cardiac-mri-1}
\end{figure}

In the remaining two figures we present a real application of the constrained level-set method to segmentation of cardiac MRI data. First we want to segment both, the left and the right ventricle of a heart. The unconstrained level-set method for the image segmentation (c.f. \cite{MikulaSarti-2007}) is not capable to separate them. In Figure~\ref{fig:cardiac-mri-1} we present a) the initial set-up of the constraint barriers; b) segmentation without constraints involved; c) segmentation with constraints separating both ventricles well. Figures~\ref{fig:cardiac-mri-1} d) and \ref{fig:cardiac-mri-1} e) demonstrate the relation between the constrained level-set method and the graph-cuts method as described in \cite{LouckyOberhuber-2012} (see the appendix for more details). Figure~\ref{fig:cardiac-mri-1} d) shows set-up of the initial seeds. They are the same as the constraints for the level-set method. There is only one more (blue) seed inside the left ventricle, otherwise it would not be segmented. The result is depicted in Figure~\ref{fig:cardiac-mri-1}.  
It is worth noting that there are other powerful methods which are specific for solving 
 cardiac segmentation problems like the one shown in Figure~\ref{fig:cardiac-mri-1}. They are based on incorporation of higher level Bayesian priors and taking into account 
 supervised learning from available cardiographic data sets (see e.g. 
 \cite{ChanZhu-2009} ) or they are based on fully unsupervised and region based level set method \cite{ParagiosDeriche-2000, Jiang-2012}

\begin{figure}
\center{
\includegraphics[width=5cm]{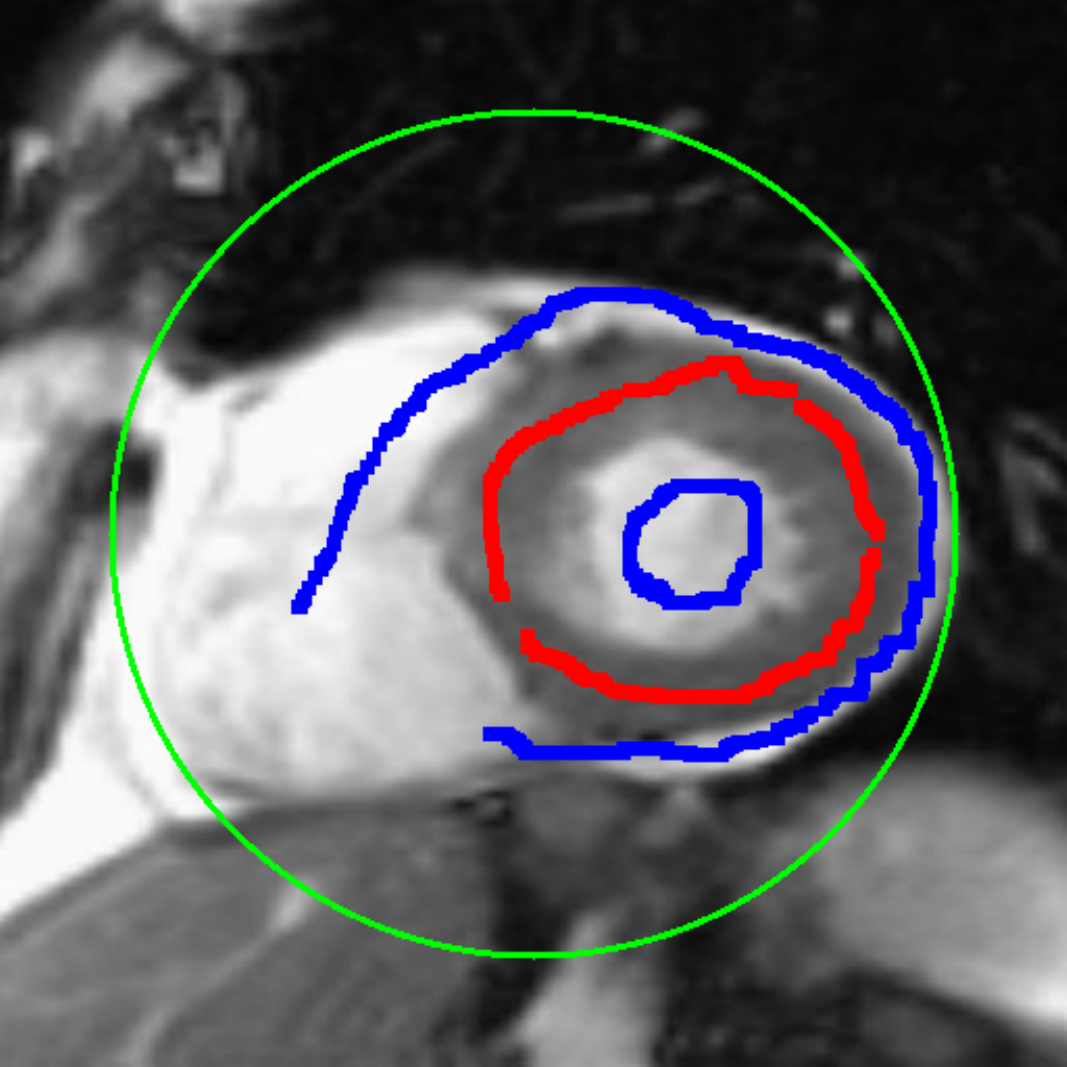}
\includegraphics[width=5cm]{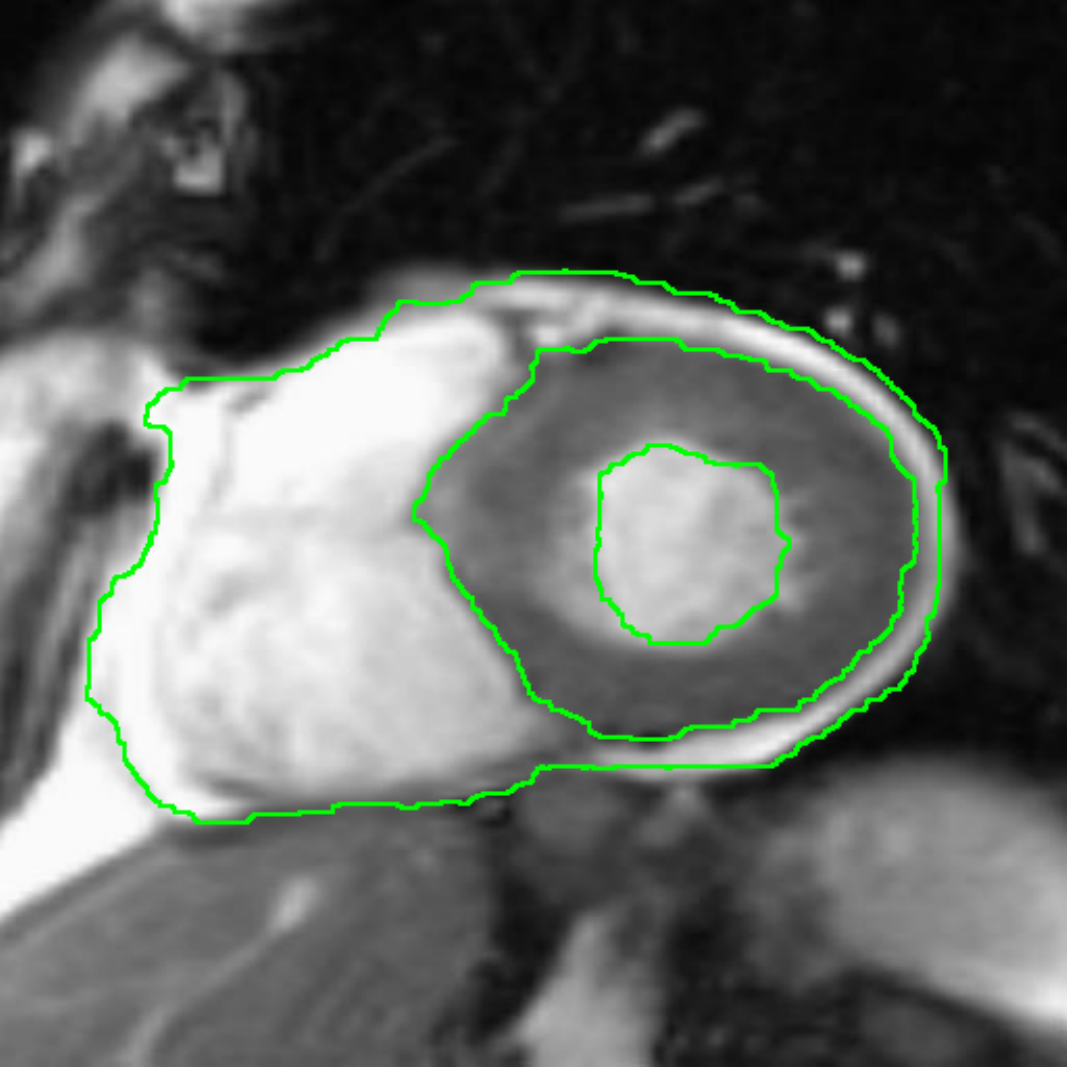}
\centerline{a) \hspace{5cm} b)} \\
\includegraphics[width=5cm]{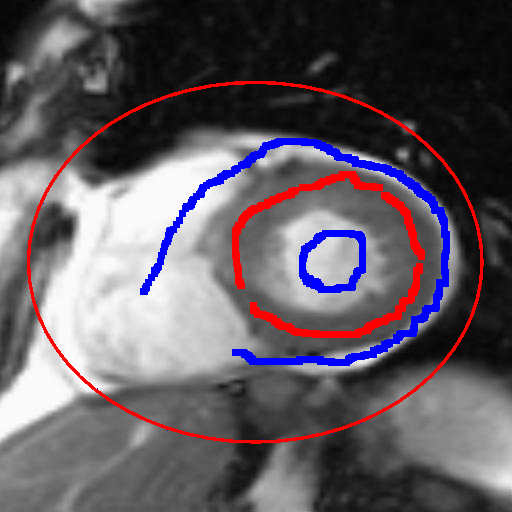}
\includegraphics[width=5cm]{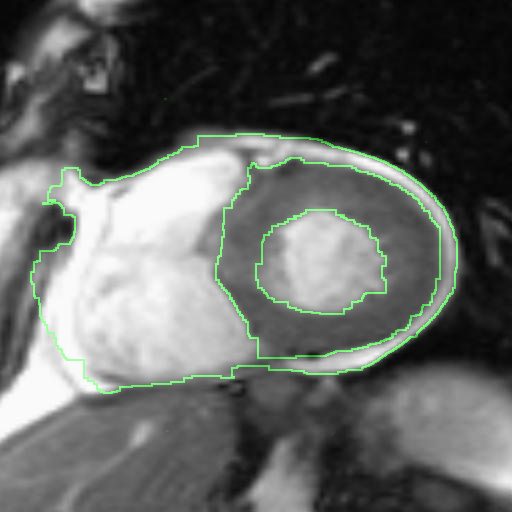}
\centerline{c) \hspace{5cm} d)} 
}
\caption{a) Setup for the constrained level-set method. b) Segmentation separating the left ventricle (small almost circular curve inside) surrounded by septum (larger curve) and the right ventricle with pericardium fat (the outer curve). c) Initial seeds for the graph-cuts method. d) Segmentation computed by the graph-cuts method.}
\label{fig:cardiac-mri-2}
\end{figure}

 More complex segmentation is studied in Figures~\ref{fig:cardiac-mri-2} a) and b). Here we have separated the left ventricle, septum and the right ventricle together with pericardium fat. Comparison with the graph-cuts method can be found in Figures~\ref{fig:cardiac-mri-2} c) and d). Note that the red ellipse in Figure~\ref{fig:cardiac-mri-2} c) serves as the initial seed as well. Compared to the initial curve green circle in Figure~\ref{fig:cardiac-mri-2} a), the red ellipse does not allow the segmentation curve to grow outside. It is not true for the initial curve for the constrained level-set method. In some applications it might be advantage. Note that even though the level-set function was negative inside the green circle (Figure \ref{fig:cardiac-mri-2} a)), the red-line constraint pushed the level-set function to positive values immediately. It is no difficulty for the numerical approximation. Thanks to this fact, the initial condition need not to be compatible with the imposed constraints.

\section*{Conclusion}

We proposed a constrained level-set model for the image segmentation. The method gives a possibility to an expert to prescribe an additional information concerning the expected shape of the segmented object. In this method, we may preset fixed constraints telling us which parts of the segmented object should be inevitably inside the segmented region and some parts must remain outside. The complementary finite-volume method was used for the numerical approximation of the level-set equation. It was combined with  the Projected Successive Over Relaxation method for solving the corresponding parabolic variational inequality problem. We demonstrated the difference between the proposed and the usual level-set method on several artificial images as well as on data from magnetic resonance imaging. 

\section*{Acknowledgments}

This work was partially supported by the grant LC06052 of NCMM Research center and Research Direction Project MSM6840770010 of the Ministry of Education of the Czech Republic, the  CVUT Student Grant Agency project SGS14/206/OHK4/3T/14 and (D\v S) VEGA grant 1/0780/15. 

\section*{Appendix}

\subsection{Time-space discretization of the level-set method}
\label{appendix-1}

In this section we show details of derivation of the full time-space discretization of  (\ref{eq:level-set-segmentation}). The numerical scheme is based on the finite volume approximation. It is a slight modification of the scheme proposed by Handlovi\v{c}ov\'a \emph{et al.} in \cite{HandlovicovaMikulaSgallari-2003}.

First we evaluate the norm of smoothed gradient of the image intensity function $s_{ij} = \left| G_\sigma \ast \grad I_0 \right|_{ij}$, where $G_\sigma(x)=(2\pi\sigma^2)^{-1} \exp(-\Vert x\Vert^2/2\sigma^2)$ (see \cite{MikulaSarti-2007} for details how to effectively calculate $s_{ij}$). Then $g^0_{ij}=1/\sqrt{1+\lambda  s^2_{ij}}$ on $\clos{\omega}_h$ and we approximate $g^0_{ij,\bar{ij}}$ on the finite volume edges as follows:
$g^0_{ij,i\pm1j} = \frac{1}{2}\left(g^0_{ij} + g^0_{i\pm1j}\right),$ \ 
$g^0_{ij,ij\pm1} = \frac{1}{2}\left(g^0_{ij} + g^0_{ij\pm1}\right)$.

Denote $\grad u^k_{ij,\bar{ij}} = \left(\partial_{x_1} u^k_{ij,\bar{ij}},\partial_{x_2} u^k_{ij,\bar{ij}}\right)$. The approximation of $\grad u^k_{ij,\bar{ij}}$ in the direction of the vector  $\nu_{ij,\bar{ij}}$ is obvious:
\begin{equation}
\partial_{x_1} u^k_{ij,i\pm1j} = \frac{u^k_{i\pm1j} - u^k_{ij}}{h_1}, \quad
\partial_{x_2} u^k_{ij,ij\pm1} = \frac{u^k_{ij\pm1} - u^k_{ij}}{h_2}.
\end{equation}
In order to calculate remaining coordinates of $\grad u^k_{ij,\bar{ij}}$ which are perpendicular to $\nu_{ij,\bar{ij}}$ we need to know the value $u^k$ at the ends of $l_{ij,\bar{ij}}$ (corners of the finite volume $v_{ij}$). They can be only approximated using the value of $u^k$ at the neighboring volumes:
\[
u^{k}_{ij,pq} = \frac{1}{4} \left( u^{k}_{ij} + u^{k}_{pj} + u^{k}_{iq} + u^{k}_{pq}\right),
\]
where $p=i\pm1, q=j\pm1$. Hence we obtain the approximation of $\grad u^k_{ij,\bar{ij}}$ in the direction perpendicular to $\nu_{ij,\bar{ij}}$ in the form 
\begin{equation}
\partial_{x_1} u^{k}_{ij,ij\pm1} = \frac{u^{k}_{ij,i+1j\pm1} - u^{k}_{ij,i-1j\pm1}}{h_1} \quad , 
\partial_{x_2} u^{k}_{ij,i\pm1j} = \frac{u^{k}_{ij,i\pm1,j+1} - u^{k}_{ij,i\pm1j-1}}{h_2}.
\end{equation}
Having calculated approximation of $\grad u^k_{ij,\bar{ij}}$ we can approximate the term $Q^k_{ij,\bar{ij}}$ as follows:
\[
Q^{k-1}_{ij,pq} = \sqrt{\epsilon^2 + \left( \partial_{x_1} u^{k-1}_{ij,pq} \right)^2 + \left( \partial_{x_2} u^{k-1}_{ij,pq} \right)^2},
\]
where $p=i\pm 1, q=j\pm1$. The "capacity" term $1/Q^k$ and the function $u^k$ are approximated by a constant value on the finite volume $v_{ij}$. Taking the averaged value 
\[
Q^{k-1}_{ij} = \frac{1}{4} \left(Q^{k-1}_{ij,i+1j}+Q^{k-1}_{ij,ij+1}+Q^{k-1}_{ij,i-1j}+Q^{k-1}_{ij,ij-1}\right),
\]
yields the approximation of the left-hand side of (\ref{eq:level-set-segmentation}):
\[
\int_{\Omega_{ij}}\frac{1}{Q^{k-1}} \frac{u^{k}-u^{k-1}}{\tau} \approx  h_1h_2\frac{1}{Q^{k-1}_{ij}} \frac{u_{ij}^{k}-u_{ij}^{k-1}}{\tau}.
\]
Applying the Stokes theorem to the right-hand side of (\ref{eq:level-set-segmentation}) yields
\begin{equation}
\label{stokes}
\int_{\Omega_{ij}}\Div{g^0 \frac{\grad u^k}{Q^{k-1}}} \dx = 
\sum_{\bar{ij} \in \mathcal{N}_{ij}}
   \int_{\Gamma_{ij,\bar{ij}}}
      \frac{g^0}{Q^{k-1}} \frac{\partial u^{k}}{\partial \nu} \dS
\approx
\sum_{\bar{ij} \in \mathcal{N}_{ij}}
l_{ij,\bar{ij}}
  \frac{g^0_{ij,\bar{ij}}}{Q^{k-1}_{ij,\bar{ij}}}
  \grad u^{k}_{ij,\bar{ij}} \cdot \nu_{ij,\bar{ij}}
, 
\end{equation}
where $\nu$ is the outer unit normal of $\Gamma_{ij}$. Here we have assumed that the fluxes are constant on each segment $\Gamma_{ij,\bar{ij}}$ of the boundary of a finite volume.
Notice that $\nu_{ij,\bar{ij}} = (\bar{i} - i, \bar{j} - j)$ (see Fig.~\ref{fig:fin-volume}). In such a regular grid, one coordinate of $\nu$ is always vanishing. It cancels one coordinate of $\grad u^k$ in the inner product $\grad u^k \cdot \nu$. Recall that $l_{ij,\bar{ij}}$ attains only the values $h_1$ or $h_2$. We have
\begin{eqnarray*}
\int_{\Omega_{ij}}\Div{g^0 \frac{\grad u^k}{Q^{k-1}}} \dx &\approx&
h_2\frac{g^0_{ij,i+1j}}{Q^{k-1}_{ij,i+1j}} \frac{u^k_{i+1j} - u^k_{ij}}{h_1} +
h_1\frac{g^0_{ij,ij+1}}{Q^{k-1}_{ij,ij+1}} \frac{u^k_{ij+1} - u^k_{ij}}{h_2}  \nonumber \\
&&+ h_2\frac{g^0_{ij,i-1j}}{Q^{k-1}_{ij,i-1j}} \frac{u^k_{i-1j} - u^k_{ij}}{h_1} +
h_1\frac{g^0_{ij,ij-1}}{Q^{k-1}_{ij,ij-1}} \frac{u^k_{ij-1} - u^k_{ij}}{h_2}.
\end{eqnarray*}
It leads to the system of linear equation (\ref{eq:algebraic-segmentation}) where 
\begin{eqnarray}
\label{def:A-1}
A_{i\pm1j}^{k} &=& -\tau Q_{ij}^{k-1} \frac{g^0_{ij,i\pm1j}}{h_1^2 Q_{ij,i\pm1j}^{k-1}}, \quad
A_{ij\pm1}^{k} = -\tau Q_{ij}^{k-1} \frac{g^0_{ij,ij\pm1}}{h_2^2 Q_{ij,ij\pm1}^{k-1}}, 
\\
\label{def:A-6}
A_{ij}^{k}   &=& 1 - \left( A_{i+1j}^{k} + A_{ij+1}^{k} + A_{i-1j}^{k} + A_{ij-1}^{k} \right).
\end{eqnarray}
Approximation of the Neumann boundary condition $\grad u \cdot \nu = 0$ on $\partial \Omega$ yields $ A^k_{0j} = A^k_{i0} = A^k_{N_1j} = A^k_{iN_2} = 1, A^k_{1j} = A^k_{i1} = A^k_{N_1-1j} = A^k_{iN_1-1} = -1$.

\subsection{Set-up of the graph-cuts method.}

\begin{table}[h]
\centering
{
\begin{tabular}{ccl@{\hspace{7mm}}c}
\toprule
\textbf{Link type} & \textbf{Edge} & & \textbf{Capacity}\\ \toprule
\emph{n}-link & $(p,q)$ & for $p,q\in \mathcal{P}$, $d=\lVert(p,q)\rVert$ & $B(\Delta I,d)$\\ \midrule
\emph{t}-link & $(s,p)$ & for $p\in \mathcal{P}\setminus\left\{\mathcal{O} \cup \mathcal{P}\right\}$ & $\lambda R_s(I_p)$\\
 & & for $p\in \mathcal{O}$ & $\infty$\\
 & & for $p\in \mathcal{P}$ & $0$\\ \cmidrule{2-4}
 & $(p,t)$ & pro $p\in \mathcal{P}\setminus\left\{\mathcal{O} \cup \mathcal{P}\right\}$ & $\lambda R_t(I_p)$\\
 & & for $p\in \mathcal{O}$ & $0$\\
 & & for $p\in \mathcal{P}$ & $\infty$\\
\bottomrule
\end{tabular}}
\caption{Summary of edge weights for the graph-cuts method. $I_p$, $I_q$ denotes intensities if pixels $p$ and  $q$ respectively; $d$ denotes a distance of the pixels.}
\label{tab:edge-weights}
\end{table}

The details of the implementation of the graph-cuts method is described in \cite{LouckyOberhuber-2012}. For the purpose of this paper we have used different set-up of the edge capacities adopted from \cite{BoykovJolly-2001}. The model is briefly described in Table (\ref{tab:edge-weights}) where

\begin{equation}\label{eq:new}
\begin{aligned}
B(\Delta I,d)&=\exp\left(-\frac{\Delta I\,{}^2}{2 \sigma_n^2}\right)\cdot\frac{1}{d}\\
R_s(I)&=-\ln P(I|\mathcal{O})\\
R_t(I)&=-\ln P(I|\mathcal{P})\\
\Delta I &=\left|I_p-I_q\right|\,.
\end{aligned}
\end{equation}

Here $\ln P(I|O)$ and $\ln P(I|P)$ are conditioned probabilities expressing if given pixel belongs to the object $O$ and background $P$ respectively. These probabilities are determined from histograms measured on the initial seeds inside and outside the object of interest respectively.

\end{document}